\documentclass[11pt,a4paper,leqno,twoside]{amsart}
\usepackage{latexsym,amssymb,amsmath}
\input xy
\xyoption{all}
\def\CC{\mathbb C}

\def\RR{\mathbb R}
\def\HH{\mathbb H}
\def\AA{{\mathbb A}}

\def\OO{\mathbb O}

\def\ZZ{\mathbb Z}

\def\NN{\mathbb N}
\def\11{\mathbf 1}
\def\PP{\mathbb P}
\def\QQ{\mathbb Q}

\def\e1{\varepsilon_1}
\def\e2{\varepsilon_2}
\def\e3{\varepsilon_3}

\def\P2{{\PP}^2}

\def\00{\underline{0}}
\def\J0{{\cal J}_3(\underline{0})}

\def\PJ0{\PP({\cal J}_3(\underline{0}))}

\def\e{\varepsilon}

\def\AP2{{\AA\PP}^2}
\def\RP2{{\RR\PP}^2}
\def\CP2{{\CC\PP}^2}
\def\HP2{{\HH\PP}^2}
\def\OP2{{\OO\PP}^2}

\newtheorem{theo}{Theorem}[section]
\newtheorem{coro}[theo]{Corollary}
\newtheorem{lemm}[theo]{Lemma}
\newtheorem{prop}[theo]{Proposition}
\theoremstyle{remark}
\newtheorem{rema}[theo]{Remark}
\theoremstyle{definition}
\newtheorem{defn}{Definition}

\begin{document}

\title{Algebraic extensions of global fields with
one-dimensional local class field theory}

\author
{I.D. Chipchakov} \email{chipchak@math.bas.bg}
\address{Institute of Mathematics and Informatics\\Bulgarian Academy
of Sciences\\Acad. G. Bonchev Str., bl. 8\\1113, Sofia, Bulgaria}
\thanks{The author was partially supported by Grant MI-1503/2005 of the
Bulgarian NSF}
\begin{abstract}
Let $E$ be an algebraic extension of a global field $E _{0}$ with a
nontrivial Brauer group Br$(E)$, and let $P(E)$ be the set of those
prime numbers $p$, for which $E$ does not equal its maximal
$p$-extension $E(p)$. This paper shows that $E$ admits
one-dimensional local class field theory if and only if there exists
a system $V(E) = \{v(p)\colon \ p \in P(E)\}$ of (nontrivial)
absolute values, such that $E(p) \otimes _{E} E _{v(p)}$ is a field,
where $E _{v(p)}$ is the completion of $E$ with respect to $v(p)$.
When this occurs, we determine by $V(E)$ the norm groups of finite
extensions of $E$, and the structure of Br$(E)$. It is also proved
that if $P$ is a nonempty set of prime numbers and $\{w(p)\colon \ p
\in P\}$ is a system of absolute values of $E _{0}$, then one can
find a field $K$ algebraic over $E _{0}$ with such a theory, so that
$P(K) = P$ and the element $\kappa (p) \in V(K)$ extends $w(p)$, for
each $p \in P$.
\end{abstract}

\maketitle

\centerline{\bf Introduction}

\par
\medskip
The purpose of this paper is to characterize the fields pointed out
in the title as well as to describe the norm group $N(R/E)$ of an
arbitrary finite extension $R$ of such a field $E$, and to determine
the structure of the Brauer group Br$(E)$. When Br$(E) \neq \{0\}$,
our research proves the validity of the Hasse principle for norm
form polynomials of $R/E$. It shows that $N(R/E) = N(\Phi (R)/E)$,
for some finite abelian extension $\Phi (R)$ of $E$ (by \cite{Ch6},
Theorem 1.2, this is not necessarily true when $E$ is an arbitrary
field admitting one-dimensional local class field theory (abbr,
LCFT)). The field $\Phi (R)$ is uniquely determined up-to an
$E$-isomorphism by the local behaviour of $R/E$ at the elements of a
certain characteristic system $\{v(p)\colon  p \in P(E)\}$ of
absolute values of $E$, indexed by the set $P(E)$ of those prime
numbers $p$, for which $E$ is properly included in its maximal
$p$-extension $E(p)$ in a separable closure $E _{\rm sep}$ of $E$
(see Theorems 2.1 and 2.3). The extension $\Phi (R)/E$ is subject to
the classical local reciprocity law, i.e. the quotient group of the
multiplicative group $E ^{\ast }$ by $N(\Phi (R)/E)$ is isomorphic
to the Galois group $\mathcal{G}(\Phi (R)/E)$. At the same time,
this paper proves that if $E _{0}$ is a global field and $P$ is a
set of prime numbers, then the class $\Omega _{P} (E _{0})$ of
algebraic extensions $\Sigma $ of $E _{0}$, which admit LCFT and
satisfy the equality $P(\Sigma ) = P$, is nonempty. It shows that
each $K \in \Omega _{P} (E _{0})$ includes as a subfield a unique
minimal element $R(K)$ of $\Omega _{P} (E _{0})$ (see Theorem 2.2
and Corollary 4.4). The established properties of the elements of
$\Omega _{P} (E _{0})$ extend the arithmetic basis of LCFT and
thereby allow one to determine the scope of validity of the norm
limitation theorem in LCFT (cf. \cite{Ch3}, Sect. 4, and \cite{Ch6},
Theorem~1.1 (ii) and Remark~3.3).
\par
\medskip
Note that a field $E$ is said to admit LCFT, if its finite abelian
extensions in $E _{\rm sep}$ are uniquely determined by their norm
groups, and for each pair $(M _{1}, M _{2})$ of such extensions, the
norm group over $E$ of the compositum $M _{1}M _{2}$ is equal to the
intersection $N(M _{1}/E) \cap N(M _{2}/E)$ and $N(M _{1} \cap M
_{2}/E) = N(M _{1}/E)N(M _{2}/E)$. Clearly, this is the case if and
only if $E$ admits local $p$-class field theory (i.e. finite abelian
extensions of $E$ in $E(p)$ have the same properties), for each $p
\in P(E)$. By \cite{Ch5}, Theorem 3.1, $E$ admits local $p$-class
field theory, for a given $p \in P(E)$, if the $p$-component Br$(E)
_{p}$ of Br$(E)$ is nontrivial, and $E$ is a $p$-quasilocal field,
i.e. its cyclic extensions of degree $p$ are embeddable as
$E$-subalgebras in each central division $E$-algebra of Schur index
$p$. By \cite{Ch5}, Remark~3.4 (ii), the converse holds when $E$ is
an algebraic extension of a global field $E _{0}$. Then $E$ admits
LCFT if and only if it is a strictly primarily quasilocal field (a
strictly PQL-field), i.e. $E$ is $p$-quasilocal with Br$(E) _{p}
\neq \{0\}$, for every $p \in P(E)$. As noted above, the main
results of the present paper are obtained by characterizing several
types of such fields in terms of valuation theory and thereafter by
a local-to-global approach.
\par
\medskip
Throughout the paper, algebras are assumed to be associative with a
unit, simple algebras are supposed to be finite-dimensional over
their centres, Brauer groups of fields are viewed as additively
written, and Galois groups are regarded as profinite with respect to
the Krull topology. The set of prime numbers is denoted by $\PP $,
and for any abelian torsion group $T$, $T _{p}$ stands for the
$p$-component of $T$, for each $p \in \PP $. Absolute values of
fields are supposed to be nontrivial, and for each algebra $A$, the
considered subalgebras of $A$ contain its unit. As usual, a field
$E$ is called formally real, if $-1$ is not presentable over $E$ as
a finite sum of squares; we say that $E$ is a nonreal field,
otherwise. In what follows, $M(E)$ denotes a system of
representatives of the equivalence classes of absolute values of
$E$, and for each $n \in \NN $, $E ^{\ast n} = \{\alpha ^{n}\colon \
\alpha \in E ^{\ast }\}$. For any field extension $\Phi /E$, $I(\Phi
/E)$, Br$(\Phi /E)$ and $\rho _{\Phi /E}$ denote the set of
intermediate fields of $\Phi /E$, the relative Brauer group of $\Phi
/E$, and the scalar extension map of Br$(E)$ into Br$(\Phi )$,
respectively. When $\Phi /E$ is finite, we write $P(\Phi /E)$ for
the set of all $p \in \PP $ dividing the degree $[\Phi \colon E]$.
The class of central simple $E$-algebras is denoted by $s(E)$,
$d(E)$ stands for the class of division algebras $\Delta \in s(E)$,
and for each $A \in s(E)$, $[A]$ is the equivalence class of $A$ in
Br$(E)$. By a group formation, we mean a nonempty class $\chi $ of
finite groups satisfying the following conditions:
\par
$\chi $ is closed under the formation of subgroups and homomorphic
images;
\par
If $G$ is a finite group, and $H _{1}$, $H _{2}$ are normal
subgroups of $G$, such that $G/H _{j} \in \chi $: $j = 1, 2$, then
$G/(H _{1} \cap H _{2}) \in \chi $.
\par \noindent
The group formation $\chi $ is called abelian closed, if it is
closed under the formation of subgroups, finite direct products and
group extensions with abelian kernels (for examples of such
formations, see, e.g., \cite{Ch6}, Remark~6.1). When this occurs,
$\chi $ is saturated, i.e. a finite group $G$ lies in $\chi $,
provided that $G/\Phi (G) \in \chi $, where $\Phi (G)$ is the
Frattini subgroup of $G$. Our terminology in valuation theory,
simple algebras, Brauer groups and abstract abelian groups is
standard (as used, e.g., in \cite{B2}; \cite{L1}; \cite{CF};
\cite{P}; \cite{M} and \cite{F}), as well as the one concerning
profinite groups, Galois cohomology, field extensions and Galois
theory (cf. \cite{S}; \cite{K}; \cite{Ko} and \cite{L1}).
\par
\medskip
The paper is organized as follows: Section 1 includes preliminaries
on algebraic extensions of a global field $E _{0}$, used in the
sequel. The main results of the paper are stated in Section 2, and
proved in Sections 3, 4 and 5. Section 3 presents a characterization
of the $p$-quasilocal fields in the class of algebraic extensions of
$E _{0}$, and a proof of the local reciprocity law for strictly
PQL-fields algebraic over $E _{0}$. In Section 4, we characterize
the minimal algebraic strictly PQL-extensions of $E _{0}$, and we
describe the Brauer groups of strictly PQL-fields algebraic over $E
_{0}$. In Section 5 we prove the Hasse norm theorem for a finite
Galois extension $M/E$, where $E \in I(\overline E _{0}/E _{0})$ and
the group $_{p} {\rm Br} (E) = \{b \in {\rm Br}(E)\colon \ pb = 0\}$
is finite, for every $p \in P(E)$. This, applied to the case where
$E$ is strictly PQL, allows us to establish the claimed properties
of norm groups of $E$.

\section {Preliminaries on the local behaviour of algebraic extensions
of global fields}
\par
\medskip
Let $E _{0}$ be a global field, $\overline E _{0}$ an algebraic
closure of $E _{0}$ and $A$ be a finite dimensional $E$-algebra, for
some $E \in I(\overline E _{0}/E _{0})$. Suppose that $B$ is a basis
of $A$, $\Sigma $ is a finite subset of $A$, $\Sigma _{1} (B)$ is
the set of structural constants of $A$, and $\Sigma _{2}(B)$ is the
one of coordinates of elements of $\Sigma $, relative to $B$. Then
the extension $E _{1}$ of $E _{0}$ generated by the union $\Sigma
_{1} (B) \cup \Sigma _{2} (B)$ is finite, and the following
statements are true (cf. [P, Sect. 9.2, Proposition c]):
\par
\medskip
(1.1) (i) The subring $A _{1}$ of $A$ generated by $E _{1} \cup B$
is an $E _{1}$-subalgebra of $A$, such that the $E$-algebras $A _{1}
\otimes _{E _{1}} E$ and $A$ are isomorphic and $\Sigma \subseteq A
_{1}$;
\par
(ii) If $A/E$ is a Galois extension and $\Sigma $ contains the roots
in $A$ of the minimal polynomial over $E _{0}$ of a given primitive
element of $A/E$, then $A _{1}/E _{1}$ is Galois with $\mathcal{G}(A
_{1}/E _{1}) \cong \mathcal{G}(A/E)$ (canonically);
\par
(iii) If $A \in d(E)$, then $\Sigma $ can be chosen so that $A _{1}
\in d(E _{1})$, exp$(A _{1}) = {\rm exp}(A)$ and ind$(A _{1}) = {\rm
ind}(A)$.
\par
\medskip
Statement (1.1) enables one to generalize some known properties of
central division algebras over global fields and local fields, as
follows:
\par
\medskip
(1.2) If $E _{0}$ is a global or locally compact field and $E \in
I(\overline E _{0}/E _{0})$, then every $D \in d(E)$ is cyclic with
ind$(D) = {\rm exp}(D)$; in addition, if $E _{0}$ is locally
compact, then $\rho _{E _{0}/E}$ is surjective.
\par
\medskip
An absolute value $v$ of a field $E$ is said to be $\chi $-Henselian
with respect to a group formation $\chi $, if $v$ is nonarchimedean
and uniquely extendable to an absolute value of the compositum
$E(\chi )$ of finite Galois extensions of $E$ in $E _{{\rm sep}}$
with Galois groups belonging to $\chi $ (or, equivalently, if $v$
extends uniquely on each finite Galois extension $M$ of $E$ with
$\mathcal{G}(M/E) \in \chi $). It follows from Galois theory and
general properties of valuation prolongations (cf. \cite{B2}, Ch.
VI) that $E$ is $\chi $-Henselian if and only if it is $\bar \chi
$-Henselian, where $\bar \chi $ is the minimal saturated group
formation including $\chi $. This allows us to restrict our
considerations of $\chi $-Henselity to the special case where $\chi
$ is saturated. When $\chi $ is the formation of $p$-groups, for a
given $p \in \PP $, the valuation $v$ is called $p$-Henselian. We
say that $v$ is Henselian, if it has a unique prolongation $v _{F}$
on each $F \in I(\overline E/E)$, i.e. if $v$ is Henselian with
respect to the formation of all finite groups. It is known
\cite{Er}, \cite{W}, that in this case, every finite dimensional
division $E$-algebra has a valuation extending $v$ (we also denote
it by $v$). The following statements describe some basic relations
between central division algebras over a Henselian valued field $(E,
v)$, and the corresponding algebras over the completion $E _{v}$
(see, for example, \cite{Ch4}):
\par
\medskip
(1.3) (i) The completion $D _{v}$ of each finite dimensional
division $E$-algebra $D$ with a centre $Z(D) \in I(E _{\rm sep}/E)$
is $E _{v}$-isomorphic to $D \otimes _{E} E _{v}$;
\par
(ii) The mapping $\rho _{E/E _{v}}$ is an isomorphism;
\par
(iii) A finite separable extension $L$ of $E$ embeds in an algebra $U \in
d(E)$ if and only if $L _{v}$ embeds in $U _{v}$ over $E _{v}$.
\par
\medskip
When $v$ is Henselian, finite extensions of $E _{v}$ in $E _{v,{\rm
sep}}$ can be characterized as follows (cf. \cite{B2}, Ch. VI, Sect.
8, No 2, Cor. 2, and \cite{Ch4}):
\par
\medskip
(1.4) Every finite separable extension $L$ of $E _{v}$ is $E
_{v}$-isomorphic to $\widetilde L _{v}$, where $\widetilde L$ is the
separable closure of $E$ in $L$. In order that $L/E _{v}$ is a
Galois extension it is necessary and sufficient that so is
$\widetilde L/E$; when this holds, $\mathcal{G}(L/E _{v}) \cong
\mathcal{G}(\widetilde L/E)$.
\par
\medskip
The presented results will usually be applied to the case of an
algebraic extension $E$ of a local field $E _{0}$. It is well-known
that the natural absolute value, say $v _{0}$, of $E _{0}$ is
Henselian, whence its unique prolongation $v$ on $E$ preserves this
property. Also, it follows from (1.1) and the classical local class
field theory (cf. \cite{CF}, Ch. VI, Sect. 1) that finite extensions
of $E$ and central division $E$-algebras are related as follows:
\par
\medskip
(1.5) (i) A finite extension $L$ of $E$ embeds in an algebra $\Delta
\in d(E)$ if and only if $[L\colon E]$ divides ind$(\Delta )$;
$[\Delta ] \in {\rm Br}(L/E)$ if and only if ind$(\Delta ) \mid
[L\colon E]$;
\par
(ii) The mapping $\rho _{E _{0}/E}$ is surjective;
\par
(iii) The fields $E$ and $E _{v}$ admit local $p$-class field
theory, for a given $p \in P(E)$, if and only if Br$(E) _{p} \neq
\{0\}$; in order that Br$(E) _{p} = \{0\}$ it is necessary and
sufficient that $E$ contains as a subfield an extension $E _{n}$ of
$E$ of degree divisible by $p ^{n}$, for each $n \in \NN $:
\par
(iv) If Br$(E) _{p} \neq \{0\}$, then there exists $K _{p} \in I(E/E
_{0})$, such that $[K _{p}\colon E _{0}] \in \NN $ and $p \dagger
[L\colon K _{p}]$, for any $L \in I(E/K _{p})$ with $[L\colon K
_{p}] \in \NN $; also, $\rho _{K _{p}/E}$ maps Br$(K _{p}) _{p}$
bijectively on Br$(E) _{p}$; in particular, Br$(E) _{p}$ is
isomorphic to the quasicyclic $p$-group $\ZZ (p ^{\infty })$, i.e.
to the $p$-component of the quotient group $\QQ /\ZZ $ of the
additive group of rational numbers by the subgroup of integers.
\par
\medskip
The following lemma is essentially a special case of the
Grunwald-Wang theorem \cite{Wa}, (see also \cite{AT}, Ch. 10,
Theorem 5, and \cite{LR}).
\par
\medskip
\begin{lemm} Assume that $E$ is an algebraic extension of a
global field $E _{0}$, $\{v _{1},..., v _{t}\}$ is a finite subset
of $M(E)$, $\widetilde F _{1},..., \widetilde F _{t}$ are cyclic
extensions of $E _{v _{1}},..., E _{v _{t}}$, respectively, of
degrees not divisible by $4$, and $n = {\rm l.c.m.}\{[\widetilde F
_{i}\colon $ $E _{v _{i}}]\colon \ i = 1,..., t\}$. Then there is a
cyclic extension $F/E$ with $[F\colon E] = n$ and $F _{v _{i}'}
\cong \widetilde F _{i}$ over $E _{v}$ whenever $i \in \{1,...,
t\}$, $v _{i} ^{\prime } \in M(F)$ and $v _{i} ^{\prime }$ extends
$v _{i}$.
\end{lemm}
\par
\medskip
\begin{proof} It follows from Krasner's lemma (cf. \cite{L2},
Ch. II, Proposition 3) that each $\tilde F _{i}/E _{v _{i}}$, $i \le
t$, has a primitive element say $\alpha _{i}$ that is algebraic over
the closure of $E _{0}$ in $E _{v _{i}}$. Hence, by (1.1) (ii),
$I(E/E _{0})$ contains a finite extension $L$ of $E _{0}$, whose
absolute values $w _{1},..., w _{t}$ induced by $v _{1},..., v
_{t}$, respectively, are pairwise nonequivalent, and for each $i \le
t$, there exists an $E _{v _{i}}$-isomorphism $\widetilde F _{i}
\cong \widetilde \Phi _{i} \otimes _{L _{w _{i}}} E _{v _{i}}$, for
some cyclic extension $\widetilde \Phi _{i}$ of $L _{w _{i}}$ in
$\widetilde F _{i}$ (with $[\widetilde \Phi _{i}\colon L _{w _{i}}]
= [\widetilde F _{i}\colon E _{v _{i}}]$). By Grunwald-Wang's
theorem, $L$ has a cyclic extension $\Phi $ of degree $n$ with $\Phi
_{w _{i}'} \cong \widetilde \Phi _{i}$ over $L _{w _{i}}$ whenever
$i \in \{1,..., t\}$, $w _{i} ^{\prime } \in M(\Phi )$ and $w _{i}
^{\prime }$ extends $v$. It is now easy to see that $\Phi \otimes
_{L} E := F$ is a field with the properties required by Lemma 1.1.
\end{proof}
\par
The generalization of the Grunwald-Wang theorem for independent
absolute values of arbitrary fields, contained in [LR], enables one
in conjunction with Galois theory to deduce the following result:
\par
\medskip
(1.6) If $E$ is a field and $v _{1}$, $v _{2}$ are $\Omega
$-Henselian absolute values of $E$, where $\Omega $ is an abelian
closed group formation containing finite $p$-groups whenever $p \in
\PP $ and there exists $\Lambda _{p} \in \Omega $ of order
divisible by $p$, then $v _{1}$ and $v _{2}$ are equivalent unless
$E(\Omega ) = E$.
\par
\medskip
Let now $E _{0}$ be a global field, $E \in I(\overline E _{0}/E
_{0})$, $v \in M(E)$, $v _{0}$ the absolute value of $E _{0}$
induced by $v$, $E _{v}$ a completion of $E$ with respect to $v$,
$\bar v$ the absolute value of $E _{v}$ continuously extending $v$,
$E _{0} ^{\prime }$ the closure of $E _{0}$ in $E _{v}$, $E ^{\prime
} = E.E _{0} ^{\prime }$, $v ^{\prime }$ the absolute value of $E
^{\prime }$ induced by $\bar v$, and $p$ a prime number. Applying
(1.4) and (1.5), one obtains the following:
\par
\medskip
(1.7) (i) $E ^{\prime }/E _{0} ^{\prime }$ is an algebraic extension
and $E _{0} ^{\prime }$ is a completion of $E _{0}$ with respect to
$v _{0}$; in particular, $E _{0} ^{\prime }$ is isomorphic to $\RR
$, $\CC $ or a local field;
\par
(ii) If $v$ is archimedean, then $E _{v}$ is isomorphic to $\RR $ or
$\CC $; hence, Br$(E ^{\prime }) = \{0\}$ unless $p = 2$ and $E _{v}
\cong \RR $;
\par
(iii) If $v$ is nonarchimedean, then $v ^{\prime }$ is Henselian;
it induces on $E _{0} ^{\prime }$ the continuous prolongation of
$v _{0}$;
\par
(iv) The absolute Galois groups $\mathcal{G} _{E'}$ and $\mathcal{G}
_{E _{v}}$ are continuously isomorphic, and $\rho _{E'/E _{v}}$ is
an isomorphism.
\par
\medskip
Next we present (slightly generalized) a part of the classical
Brauer-Hasse-Noether-Albert theorem (see \cite{We}, Ch. XIII, Sects.
3 and 6).
\par
\medskip
\begin{prop} Let $E _{0}$ be a global field and $E$ algebraic
extension of $E _{0}$. Then $\rho _{E/E _{v}}$ is surjective, for
each $v \in M(E)$, and the homomorphism of {\rm Br}$(E)$ into the
direct sum $\oplus _{v \in M(E)}$ Br$(E _{v})$ mapping, for each $A
\in s(E)$, $[A]$ into the sequence $[A \otimes _{E} E _{v}]\colon \
v \in M(E)$, is injective.
\end{prop}
\par
\medskip
\begin{proof} Our latter assertion is a special case of \cite{Ef},
Lemma~3.3. For the proof of the former one, fix an absolute value $v
\in M(E)$ as well as an algebra $\Delta \in d(E _{v})$, and denote
by $E _{0} ^{\prime }$ and $E ^{\prime }$ the fields from $I(E
_{v}/E _{0})$ defined in (1.7). By (1.7) (iv), there exists an $E
_{v}$-isomorphism $\Delta \cong \widetilde \Delta \otimes _{E'} E
_{v}$, for some $\widetilde \Delta \in d(E ^{\prime })$, uniquely
determined up-to an $E ^{\prime }$-isomorphism, and by (1.2),
$[\widetilde \Delta ] = [\widetilde \Delta _{0} \otimes _{E _{0}'} E
^{\prime }]$, in Br$(E ^{\prime })$, for a suitably chosen
$\widetilde \Delta \in d(E _{0} ^{\prime })$. Using consecutively
(1.7) (i) and the description of Br$(E _{0})$ by class field theory
(see, e.g., \cite{We}, Ch. XIII, Sect. 6), one obtains the existence
of an algebra $\Delta _{0} \in d(E _{0})$, such that $[\Delta _{0}
\otimes _{E _{0}} E _{0} ^{\prime }] = [\widetilde \Delta _{0}]$ (i
Br$(E _{0} ^{\prime })$). These observations and the general
behaviour of tensor products under scalar extensions (cf. \cite{P},
Sect. 9.4, Corollary a) indicate that $[\Delta ] = [(\Delta _{0}
\otimes _{E _{0}} E) \otimes _{E} E _{v}]$, in Br$(E _{v})$, which
proves the former assertion of Proposition 1.2.
\end{proof}
\par
\medskip
Applying (1.7) (ii), Proposition 1.2 and the latter part of (1.5)
(iii), one proves the following result of Fein and Schacher
\cite{FS}, Sect. 2, Theorem 4):
\par
\medskip
(1.8) If $E _{0}$ is a global or local field, $E \in I(\overline E
_{0}/E _{0})$ and Br$(E) _{p} = \{0\}$, for some $p \in \PP $, then
Br$(E _{1}) _{p} = \{0\}$, for every $E _{1} \in I(\overline E
_{0}/E)$.
\par
\medskip
\begin{lemm} Let $F$ be a field, $w \in M(F)$, $M/F$ a
Galois extension, and $P _{w} (M)$ the set of prolongations of $w$
on $M$. Then the action of $\mathcal{G}(M/F)$ on $P _{w} (M)$ by the
rule $v \to v \circ \varphi \colon \ \varphi \in \mathcal{G}(M/F)$,
$v \in P _{w} (M)$, is transitive.
\end{lemm}
\par
\medskip
\begin{proof} The assertion is well-known in case
$[M\colon F] \in \NN $ (cf. \cite{L1}, Ch. IX, Proposition 11), so
we assume that $[M\colon F] = \infty $. Fix prolongations $v _{1}$
and $v _{2}$ of $w$ on $M$, denote by $\Sigma (M/F)$ the set of
finite Galois extensions of $F$ in $M$, and for each $L \in \Sigma
(M/F)$, let $v _{1} (L)$ and $v _{2} (L)$ be the absolute values of
$L$ induced by $v _{1}$ and $v _{2}$, respectively. We show that the
sets $P(L) = \{\sigma \in \mathcal{G}(M/F)\colon \ \sigma (v _{1})$
extends $v _{2} (L)$ are closed, for all $L \in \Sigma (M/F)$, and
their intersection $P$ is nonempty. It follows from Galois theory
that $P(L)$ is a coset of a subgroup $U _{L}$ of $\mathcal{G}(M/F)$
including $\mathcal{G}(M/L)$, for each $L \in \Sigma (M/F)$.
Observing also that $P(L _{1}... L _{s}) \subseteq \cap _{j=1} ^{s}
P(L _{j})$ whenever $s \in \NN $ and $L _{1},..., L _{s} \in \Sigma
(M/F)$, one obtains from the compactness of $\mathcal{G}(M/F)$ that
$P \neq \phi $. Since $g(v _{1}) = v _{2}$, for every $g \in P$,
this proves Lemma 1.3.
\end{proof}

\section{Statements of the main results}

\par
\medskip
It is clear from Galois theory and the general description of
relative Brauer groups of cyclic extensions [P, Sect. 15.1,
Proposition b] that a field $E$ with Br$(E) = \{0\}$ is strictly PQL
if and only if $P(E) = \phi $, i.e. $E$ does not possess proper
abelian extensions. The following result characterizes the strictly
PQL-fields $E$ with $P(E) \neq \phi $ in the set $I(\overline E
_{0}/E _{0})$, for a global field $E _{0}$:
\par
\medskip
\begin{theo} Let $E _{0}$ be a global field and $E$ an extension of
$E _{0}$ in $\overline E _{0}$, such that $P(E) \neq \phi $. Then
the following conditions are equivalent:
\par
{\rm (i)} $E$ is a strictly PQL-field;
\par
{\rm (ii)} For each $p \in P(E)$, {\rm Br}$(E) _{p} \neq \{0\}$ and
there exists an absolute value $v(p)$ of $E$, such that the tensor
product $E(p) \otimes _{E} E _{v(p)}$ is a field.
\par \noindent
When these conditions hold, $v(p)$ is uniquely determined, up-to an
equivalence, for each $p \in P(E)$. Moreover, $\rho _{E/E _{v(p)}}$
maps {\rm Br}$(E) _{p}$ bijectively on {\rm Br}$(E _{v(p)}) _{p}$,
and $E(p) \otimes _{E} E _{v(p)} \cong E _{v(p)}(p)$ as $E
_{v(p)}$-algebras.
\end{theo}
\par
\medskip
\begin{defn} Under the hypotheses of Theorem 2.1, a system of absolute
values of $E$ is called characteristic, if it is subject to the
restrictions in (ii).
\end{defn}
\par
\medskip
Theorem 2.1, Proposition 1.2 and statements (1.7) (iv) and (1.8)
indicate that if $E$ is an algebraic extension of a global field $E
_{0}$, $F \in I(E/E _{0})$, and for each $p \in P(E)$, $w(p)$ is the
absolute value of $F$ induced by $v(p)$, then the groups Br$(F)
_{p}$ and Br$(F _{w(p)}) _{p}$ are nontrivial. Therefore, our
research concentrates on the study of the following class of fields:
\par
\medskip
\begin{defn} Let $E _{0}$ be a global field, $F$ an extension of $E
_{0}$ in $\overline E _{0}$ with Br$(F) \neq \{0\}$, $P$ a nonempty
set of prime numbers $p$, for which Br$(F) _{p} \neq \{0\}$ and
$\{w(p)\colon \ p \in P\}$ a system of absolute values of $F$, such
that Br$(F _{w(p)}) _{p} \neq \{0\}$. Denote by $\Omega (F, P, W)$
the set of all $E \in I(\overline E _{0}/F)$ with the following
properties:
\par
(i) $E$ is a strictly PQL-field and $P(E) = P$;
\par
(ii) $E$ possesses a characteristic system $\{v(p)\colon \ p \in P\}$
(of absolute values), such that $v(p)$ is a prolongation of $w(p)$,
for each $p \in P$.
\end{defn}
\par
\medskip
Our main results about $\Omega (F, P, W)$ can be stated as follows:
\par
\medskip
\begin{theo} With assumptions and notations being as
above, $\Omega (F, P, W)$ is a nonempty set, for which the following
assertions hold true:
\par
{\rm (i)} Every $E \in \Omega (F, P, W)$ possesses a unique subfield
$R(E)$ that is a minimal element of $\Omega (F, P, W)$;
\par
{\rm (ii)} The Galois closure over $F$ of any minimal element $E$ of
$\Omega (F, P, W)$ has a prosolvable Galois group; furthermore, $E$
is presentable as a union $\cup _{n=1} ^{\infty } K _{n}$, where $K
_{1} = E$ and $K _{n+1}$ is an intermediate field of the maximal
extension of $K _{n}$ in $E _{\rm sep}$ with a pronilpotent Galois
group;
\par
{\rm (iii)} If $E$ is a minimal element of $\Omega (F, P, W)$, $p
\in P$ and $F _{w(p)}$ is the closure of $F$ in $E _{v(p)}$, then
the degrees of finite extensions of $F _{w(p)}$ in $E _{v(p)}$ are
not divisible by $p$;
\par
{\rm (iv)} For each linearly ordered subset $\Lambda $ of $\Omega (F, P,
W)$, the intersection $\cap _{L \in \Lambda } L$ lies in $\Omega (F,
P, W)$.
\end{theo}
\par
\medskip
Let $E _{0}$ be a global field, $E \in I(\overline E _{0}/E _{0})$
and $R$ a finite extension of $E$. An element $\alpha \in E ^{\ast
}$ is called a local norm of $R/E$, if $\alpha \in N(R _{v'}/E
_{v})$, for every absolute value $v$ of $E$, and each prolongation
$v ^{\prime }$ of $v$ on $R$. Denote by $N _{\rm loc}(R/E)$ the
multiplicative group of local norms of $R/E$. As it turns out, when
$E$ is a strictly PQL-field, $N _{\rm loc} (R/E)$ is a subgroup of
$N(R/E)$ and both are determined by the local behaviour of $R/E$.
\par
\medskip
\begin{theo} Let $E _{0}$ be a global field, $E \in I(\overline E
_{0}/E _{0})$ a field admitting LCFT with $P(E) \neq \phi $, and $V
= \{v(p)\colon \ p \in P(E)\}$ a characteristic system of $E$. Also,
let $R$ be a finite separable extension of $E$ in $\overline E
_{0}$, $M$ a normal closure of $R$ in $\overline E _{0}$ over $E$,
and $\Pi (M/E) = P(M/E) \cap P(E)$. Then there exist abelian finite
extensions $F _{1}$ and $F _{2}$ of $E$ in $\overline E _{0}$, and a
real number $\varepsilon > 0$, for which:
\par
{\rm (i)} $N(F _{1}/E) = N(R/E)$, $N(F _{2}/E) = N _{\rm loc}(R/E)$,
$F _{1} \subseteq F _{2}$ and $[F _{1}\colon E]$ is not divisible by
any prime number out of $\Pi (M/E)$; the groups $E ^{\ast}/N(R/E)$
and $E ^{\ast }/N _{\rm loc}(R/E)$ are isomorphic to $\mathcal{G}(F
_{1}/E)$ and $\mathcal{G}(F _{2}/E)$, respectively;
\par
{\rm (ii)} $F _{1} = F _{2}$, provided that $\mathcal{G}(M/E)$ is
nilpotent or $\Pi (M/E) = \phi $; in the latter case, $F _{1} = F
_{2} = E$;
\par
{\rm (iii)} If $\Pi (M/E) \neq \phi $, $p \in \Pi (M/E)$ and $\Sigma
_{p}$ is the set of prolongations of $v(p)$ on $R$, then the maximal
$p$-extension of $E$ in $F _{1}$ equals the intersection $\cap
_{v(p)' \in \Sigma _{p}} A (v(p) ^{\prime })$, where $A (v(p)
^{\prime })$ is the abelian $p$-extension of $E$ in $\overline E
_{0}$ with the property that $A (v(p) ^{\prime }) \otimes _{E} E
_{v(p)}$ is isomorphic as an algebra over $E _{v(p)}$ to the maximal
abelian $p$-extension of $E _{v(p)}$ in $R _{v(p)'}$;
\par
{\rm (iv)} With assumptions being as in (iii), the maximal $p$-extension
of $E$ in $F _{2}$ equals the compositum of the fields $A (v(p)
^{\prime })\colon  v(p) ^{\prime } \in \Sigma _{p}$.
\end{theo}
\par
\medskip
Theorem 2.3 and our next result are proved at the end of Section
5.
\par
\medskip
\begin{prop} With assumptions being as in Theorem 2.3, let $f$ be a norm
form of the extension $R/E$, and $c$ an element of $E ^{\ast }$.
Then the equation $f = c$ has a solution over $E$ if and only if it
is solvable over $E _{v(p)}$, for every $p \in \Pi (M/E)$; in
particular, this is true, if $v(p) (c - 1) < \varepsilon $, for each
$p \in \Pi (M/E)$, provided that $\varepsilon > 0$ is a sufficiently
small number.
\end{prop}

\par
\medskip
\begin{rema} It is well-known that if $E$ is an algebraic
extension of a global field of characteristic $q > 0$, and $E ^{q} =
\{e ^{q}\colon \ e \in E\}$, then $[E\colon E ^{q}]$ is equal to $1$
or $q$. This implies that if $R/E$ is a finite extension and $R
_{0}$ is the separable closure of $E$ in $R$, then $N(R/R _{0}) = R
_{0} ^{\ast }$ and $N(R/E) = N(R _{0}/E)$, i.e. $N(R/E)$ can be
fully described by applying Theorem 2.3. One also sees that norm
forms of $R/E$ are subject to Hasse's principle.
\end{rema}
\par
\medskip
The following statement shows the validity the Hasse norm principle
for finite Galois extensions of strictly PQL-fields.
\par
\medskip
\begin{prop} Assume that $E _{0}$ is a global field,
$E$ is an algebraic strictly PQL-extension of $E _{0}$, $P(E) \neq
\{\phi \}$ and $V = \{v(p)\colon \ p \in P(E)\}$ is a characteristic
system of $E$. Also, let $M$ be a finite Galois extension of $E$ in
$\overline E _{0}$, and let $\Pi (M/E) = P(M/E) \cap P(E)$. Then
there exists a finite abelian extension $\widetilde M$ of $E$
satisfying the following conditions:
\par
{\rm (i)} $N(\widetilde M/E) = N(M/E) = N _{\rm loc} (M/E)$;
\par
{\rm (ii)} $[\widetilde M\colon E] \mid [M\colon E]$; in particular,
$\widetilde M = E$ in case $\Pi (M/E) = \phi $;
\par
{\rm (iii)} For each $p \in P(\widetilde M/E)$ and any absolute
value $v(p) ^{\prime }$ of $M$ extending $v(p)$, the maximal abelian
$p$-extension of $E _{v(p)}$ in $M _{v(p)'}$ is $E
_{v(p)}$-isomorphic to $\widetilde M _{p} \otimes _{E} E _{v(p)}$,
where $\widetilde M _{p} = \widetilde M \cap E(p)$.
\par
The field $\widetilde M$ is uniquely determined by $M$, up-to an
$E$-isomorphism.
\end{prop}

\par
\begin{proof} This follows from Theorem 2.3, Proposition 2.4
and Lemma 1.3.
\end{proof}

\par
\medskip
Under the hypotheses of Proposition 2.6, one obtains from Theorem
2.1 and Galois theory that $M _{0} \subseteq \widetilde M$, where $M
_{0}$ is the maximal abelian extension of $E$ in $M$. In addition,
it becomes clear that $\widetilde M = M _{0}$, provided that
$\mathcal{G}(M/E)$ is nilpotent. This is not necessarily true when
$\mathcal{G}(M/E)\cong G$, for any given nonnilpotent finite group
$G$ (cf. \cite{Ch3}, Sect. 4). Also, Proposition 2.6 shows that the
index $\vert E ^{\ast }\colon N(M/E)\vert $ divides the order $o(G)$
of $G$ and is divisible by $\vert G\colon [G, G]\vert $, $[G, G]$
being the commutator subgroup of $G$. The concluding result of this
Section characterizes those $G$, for which $M$ can be chosen so that
$\mathcal{G}(M/E) \cong G$ and $\vert E ^{\ast }\colon N(M/E)\vert =
o(G)$. It solves for such $G$ the inverse problem concerning the
admissible values of $\vert E ^{\ast }\colon N(M/E)\vert $.
\par
\medskip
\begin{prop} In the setting of Proposition 2.6, if $[M\colon E] =
[\widetilde M\colon E]$, then the Sylow $p$-subgroups of
$\mathcal{G}(M/E)$ are abelian, for every $p \in \PP $.
\par
Conversely, let $G$ be a finite group with $o(G) = m$ and abelian
Sylow subgroups, and let $n \in \NN $ divide $m$ and be divisible by
$\vert G\colon [G, G]\vert $. Then there exist subfields $E _{n}$
and $M _{n}$ of $\overline {\QQ }$, such that $E _{n}$ is a strictly
PQL-field, $M _{n}/E _{n}$ is a Galois extension with $\mathcal{G}(M
_{n}/E _{n}) \cong G$, and $\vert E _{n} ^{\ast }\colon N(M _{n}/E
_{n})\vert = n$.
\end{prop}

\begin{proof} Our former assertion can be deduced from Proposition
2.6 and Galois theory. The proof of the latter one relies on the
existence of subfields $E _{m}$ and $M _{m}$ of $\overline {\QQ }$
satisfying the following conditions (see \cite{Ch3}, Sect. 4):
\par
\medskip
(2.1) (i) $E _{m}$ is a strictly PQL-field, $P(E _{m}) := P$ equals
the set of all prime numbers, and the characteristic system $\{v
_{m} (p)\colon \ p \in P\}$ has the property that $v _{m} (p)$
induces on $\QQ $ the natural $p$-adic absolute value, for each $p
\in P$;
\par
(ii) $M _{m}/E _{m}$ is a Galois extension with $\mathcal{G}(M
_{m}/E _{m}) \cong G$, and for each $p \in P$, $M _{m,v _{m} (p)'}/E
_{m, v _{m} (p)}$ is Galois with $\mathcal{G}(M _{m,v _{m} (p)'}/E
_{m,v _{m} (p)})$ isomorphic to a Sylow $p$-subgroup of $G$, where
$v _{m} (p) ^{\prime }$ is a prolongation of $v _{m} (p)$ on $M$.
\par
\medskip
Let $P(G)$ be the set of prime divisors of $m$, and let $G _{p}$ be
a Sylow $p$-subgroup of $G$, for each $p \in P(G)$. It is clear from
Proposition 2.6, (2.1) (iii) and the commutativity of the groups $G
_{p}\colon \ p \in P(G)$, that $E _{m}$ and $M _{m}$ have the
properties required by Proposition 2.7. In particular, $[\widetilde
M _{m}\colon E _{m}] = m$, where $\widetilde M _{m}$ is the abelian
extension of $E _{m}$ associated with $M _{m}/E _{m}$, by
Proposition 2.6. We show that $\widetilde M _{m}$ possesses a
subfield $\Phi _{m}$, such that $\Phi _{m} \cap M _{m,ab} = E _{m}$
and $\Phi _{m}M _{m,ab} = \widetilde M _{m}$. Identifying as we can
(cf. \cite{K}, Ch. 7, Corollary 5.3) the character groups $C(G)$ and
$C(G _{p})$ with the cohomology groups $H ^{1} (G, \QQ /\ZZ )$ and
$H ^{1} (G _{p}, \QQ /\ZZ $, respectively, one obtains from
\cite{S}, Ch. I, Proposition 9) that the composition cor$_{p} \circ
{\rm res}_{p}$ of the restriction map res$_{p}\colon C(G) \to C(G
_{p})$ and the corestriction cor$_{p}\colon C(G _{p}) \to C(G)$
induces an automorphism of the $p$-component $C(G) _{p}$ of $C(G)$.
This implies res$_{p}$ maps $C(G) _{p}$ bijectively on a pure
subgroup of $C(G _{p})$. As pure subgroups of $C(G _{p})$ are direct
summands in $C(G _{p})$ (cf. \cite{F}, Theorem 24.5), the existence
of $\Phi _{m}$ follows now from Galois theory. It is not difficult
to see that $\Phi _{m}$ contains as a subfield an abelian extension
$F _{n}$ of $E _{m}$ of degree $m/n$, for each $n \in \NN $ subject
to the restrictions of Proposition 2.7. Observing that $F _{m} \cap
M _{m} = E _{m}$, one obtains from Galois theory that $(M _{m}F
_{m})/F _{n}$ is a Galois extension with $\mathcal{G}((M _{m}F
_{n})/F _{n}) \cong G$. Fix an absolute value $w _{n} (p)$ of $F
_{n}$ extending $v _{m} (p)$, for each $p \in P$, and put $W _{n} =
\{w _{n} (p)\colon \ p \in P\}$. Since, by (2.1) (i) and Theorem
2.1, Br$(E _{n}) \cong \QQ /\ZZ $, we have Br$(F _{n,w _{n} (p)})
_{p} \neq \{0\}$, $p \in P$, i.e. $\Omega (F _{n}, P, W _{n})$ is
well-defined. Applying Theorem 2.2 (iii), one concludes that if $M
_{n} = M _{m}E _{n}$, for some minimal element $E _{n}$ of $\Omega
(F _{n}, P, W _{n})$, then $M _{n}/E _{n}$ has the property required
by Proposition 2.7.
\end{proof}

\section{Algebraic $p$-quasilocal extensions of global fields and
Galois groups of their maximal $p$-extensions}
\par
\medskip
In this Section we characterize algebraic extensions of global
fields admitting local $p$-class field theory, for a given prime
number $p$, and thereby prove Theorem 2.1. Our argument is
presented by the following three lemmas.
\par
\medskip
\begin{lemm} Let $E$ be an algebraic extension of a global
field $E _{0}$, such that $E(p) \otimes _{E} E _{v _{1}}$ is a
field, for some $p \in \PP $ and $v _{1} \in M(E)$, and let $v _{2}
\in M(E)$ be nonequivalent to $v _{1}$. Then:
\par
{\rm (i)} $E _{v _{2}} (p) = E _{v _{2}}$ and $E(p) \otimes _{E} E _{v
_{1}}$ is $E _{v _{1}}$-isomorphic to $E _{v _{1}} (p)$;
\par
{\rm (ii)} Br$(E _{v _{2}}) _{p} = \{0\}$ and $\rho _{E/E _{v
_{1}}}$ maps {\rm Br}$(E) _{p}$ bijectively on {\rm Br}$(E _{v
_{1}}) _{p}$.
\end{lemm}
\par
\medskip
\begin{proof} It is clear from Galois theory (cf. \cite{Ko},
Proposition 2.11) that $E(p)$ embeds in $E _{v _{1}} (p)$ as an
$E$-subalgebra. Identifying $E(p)$ with its $E$-isomorphic copy in
$E _{v _{1}} (p)$, one gets from the condition on $v _{1}$ that
$E(p) \cap E _{v _{1}} = E$. Note also that every extension $L$ of
$E _{v _{1}}$ in $E _{v _{1}} (p)$ of degree $p ^{n}$ equals
$\widetilde LE _{v _{1}}$, for some $\widetilde L \in I(E(p)/E)$
with $[\widetilde L\colon E] = p ^{n}$. If $n = 1$, this is implied
by Lemma 1.1, so we assume that $n \ge 2$. It follows from the
subnormality of proper subgroups of finite $p$-groups (cf.
\cite{L1}, Ch. I, Sect. 6) and Galois theory that $L$ contains as a
subfield a cyclic extension $L _{0}$ of $E _{v _{1}}$ of degree $p$.
As $L _{0} (p) = E(p)$ and $E(p) \otimes _{L _{0}} L _{0,v'}$ is a
field, $v ^{\prime }$ being a prolongation of $v$ on $L _{0}$
(unique when $v$ is nonarchimedean), this enables one to complete
the proof of our assertion, arguing by induction on $n$. The
obtained result shows that $E _{v _{1}} (p) = E(p).E _{v _{1}}$ and
there is an $E _{v _{1}}$-isomorphism $E _{v _{1}} (p) \cong E(p)
\otimes _{E} E _{v _{1}}$. Using again Galois theory, one also
concludes that $\mathcal{G}(E(p)/E) \cong \mathcal{G}(E _{v _{1}}
(p)/E _{v _{1}})$. It is now clear from Lemma 1.1 and the condition
on $E(p)$ and $E _{v _{1}}$ that $p \notin P(E _{v _{2}})$. Hence,
by (1.2), Br$(E _{v _{2}}) _{p} = \{0\}$, so Lemma 3.1 (ii) reduces
to a consequence of Proposition 1.2.
\end{proof}

\par
\medskip
\begin{lemm} An algebraic extension $E$ of a global field
$E _{0}$ admits local $p$-class field theory, for a given $p \in
P(E)$, if and only if {\rm Br}$(E) _{p} \neq \{0\}$ and there exists
an absolute value $v$ of $E$ such that $E(p) \otimes _{E} E _{v}$
is a field. When this occurs, $v$ is unique, up-to an equivalence.
\end{lemm}

\par
\medskip
\begin{proof} It is clear from Lemma 3.1 that Br$(E) _{p} = \{0\}$
in case $E$ has nonequivalent absolute values $v _{1}$ and $v _{2}$
for which $E(p) \otimes _{E} E _{v _{1}}$ and $E(p) \otimes _{E} E
_{v _{2}}$ are fields. This proves the concluding statement of Lemma
3.2. Suppose now that $E$ admits local $p$-class field theory. This
means that $E$ is $p$-quasilocal with Br$(E) _{p} \neq \{0\}$.
Hence, by Merkurjev's theorem \cite{M}, Sect. 4, Theorem 2, there
exists $D \in d(E)$ of index $p$, by Proposition 1.2, $D \otimes
_{E} E _{v} \in d(E _{v})$, for some absolute value $v$ of $E$. We
show that $E(p) \otimes _{E} E _{v}$ is a field by assuming the
opposite. In view of Galois theory and the subnormality of proper
subgroups of finite $p$-groups, then $E _{v}$ contains as a subfield
a cyclic extension $L$ of $E$ of degree $p$. Since the algebras $D
\otimes _{E} E _{v}$ and $(D \otimes _{E} L) \otimes _{L} E _{v}$
are isomorphic over $E _{v}$ (cf. \cite{P}, Sect. 9.4, Corollary a),
this leads to the conclusion that $[D] \notin {\rm Br}(L/E)$. The
obtained result contradicts the assumption that $E$ is
$p$-quasilocal, and thereby, proves that $E(p) \otimes _{E} E _{v}$
is a field.
\par
Conversely, assume that Br$(E) _{p} \neq \{0\}$ and $E(p) \otimes
_{E} E _{v}$ is a field, for some absolute value $v$ of $E$. By
(1.2), then we have $E(p) \neq E$. Let $\overline E _{v}$ be an
algebraic closure of $E _{v}$, $M$ a cyclic extension of $E$ in
$\overline E _{v}$ of degree $p$, and $E ^{\prime }$ the
intermediate field of $E _{v}/E$ defined as in (1.7). We show that
$[\Delta ] \in {\rm Br}(M/E)$, for each $\Delta \in d(E)$ of index
$p$. It follows from Galois theory and the assumptions on $v$ and
$M$ that $ME _{v}/E _{v}$ and $ME ^{\prime }/E ^{\prime }$ are
cyclic field extensions of degree $p$. This means that $v$ is
uniquely extendable to an absolute value $v _{M}$ of $M$, and $ME
_{v}$ is a completion of $M$ with respect to $v _{M}$ (cf.
\cite{CF}, Ch. II, Theorem 10.2). Applying now (1.2), Lemma 3.1 and
Proposition 1.2, one concludes that our assertion will be proved, if
we show that $[\Delta \otimes _{E} M] \in {\rm Br}(ME _{v}/M)$.
Since the algebras $(\Delta \otimes _{E} M) \otimes _{M} (ME _{v})$,
$\Delta \otimes _{E} (ME _{v})$ and $(\Delta \otimes _{E} (ME
^{\prime })) \otimes _{ME'} (ME _{v})$ are isomorphic over $ME
_{v}$, it is sufficient to establish that $[\Delta ] \in {\rm Br}(ME
^{\prime }/E)$. This property of $\Delta $ follows from the fact
that $E ^{\prime }$ is quasilocal, the equality $[(ME ^{\prime
})\colon E ^{\prime }] = p$ and the existence of an $(ME ^{\prime
})$-isomorphism $\Delta \otimes _{E} (ME ^{\prime }) \cong (\Delta
\otimes _{E} E ^{\prime }) \otimes _{E'} (ME ^{\prime })$, so the
proof of Lemma 3.2 is complete.
\end{proof}

\par
\medskip
\begin{lemm} Let $E _{0}$ be a global field, $E$ an
algebraic extension of $E _{0}$ admitting local $p$-class field
theory, for some $p \in P(E)$, and let $v \in M(E)$ be chosen so
that $E(p) \otimes _{E} E _{v}$ is a field. Then:
\par
{\rm (i)} Br$(E) _{p}$ is isomorphic to the group $\ZZ  (p ^{\infty })$
unless $p = 2$ and $E$ is a formally real field; when this is the
case, $v$ is nonarchimedean;
\par
{\rm (ii)} Br$(E) _{2}$ is of order $2$ in case $E _{v}$ is isomorphic to
$\RR $.
\end{lemm}

\par
\medskip
\begin{proof} By \cite{Ch2}, I, Lemma 3.6, a $2$-quasilocal field $E$
is formally real if and only if Br$(E) _{2}$ is of order $2$. This,
combined with (1.5) (iv), (1.7) (i), (iv) and Lemma 3.1, proves our
assertion.
\end{proof}

\par
\medskip
Our next result shows that an algebraic extension $E$ of a global
field $E _{0}$ satisfying the condition Br$(E) \neq \{0\}$ is a
field with local class field theory in the sense of \cite{NP} if and
only if finite extensions of $E$ are strictly PQL-fields.
\par
\medskip
\begin{coro} Let $E$ be an algebraic extension of a
global field $E _{0}$ and $\Pi (E)$ the set of those prime numbers
$p$ for which $\mathcal{G} _{E}$ is of nonzero cohomological
$p$-dimension. Then finite extensions of $E$ are strictly {\rm
PQL}-fields if and only if the following assertions hold true:
\par
{\rm (i)} $E _{\rm sep} \otimes _{E} E _{v}$ is a field, for some $v
\in M(E)$;
\par
{\rm (ii)} Br$(E) _{p} \neq \{0\}$, for every $p \in \Pi (E)$.
\end{coro}

\par
\medskip
\begin{proof} Statements (1.8) and the concluding observations in
\cite{Ch5}, Sect. 1, imply finite extensions of $E$ are strictly
PQL-fields and only if $E$ is quasilocal and cd$_{p} (\mathcal{G}
_{E}) \neq 1$, for any $p \in \PP $. In view of \cite{Ch1},
Proposition 3.1, this means that if $E$ is formally real, then the
Sylow pro-$p$-subgroups of $\mathcal{G} _{E}$ are isomorphic to $\ZZ
_{p} ^{2}$, for each $p \in \Pi (E) \setminus \{2\}$. The obtained
result in fact proves that $E$ is real closed, since, by \cite{Ge},
Theorem 2.3, the algebraicity of $E/\QQ $ guarantees that
abelian subgroups of $\mathcal{G} _{E}$ are procyclic. Our argument
also relies on the fact \cite{Ch2}, I, Lemma 3.6, that $E$ is
formally real and $2$-quasilocal if and only if $[E(2)\colon E] =
2$, and when this occurs, Br$(E) _{2}$ is of order $2$. Hence, by
the Artin-Schreier theory (cf. \cite{L1}, Ch. XI, Sect. 2) and Lemma
3.3, $E$ is a real closed field if and only if it has an absolute
value $v$, such that $E _{\rm sep} \otimes _{E} E _{v}$ is a field
and $E _{v} \cong \RR $. These observations prove Corollary 3.4 in
the special case where $E$ is formally real. Assuming further that
$E$ is nonreal, we complete the proof in the following two steps:
\par
Step 1. Suppose that $E$ satisfies conditions (i) and (ii), and let
$F$ be a finite extension of $E$ in $\overline E$. Then it follows
from Lemmas 3.2 and 3.3 that $E$ is strictly PQL and Br$(E) _{p}
\cong \ZZ (p ^{\infty })$, for each $p \in \Pi (E)$. Condition (i)
is equivalent to the one that $v$ is Henselian. The prolongation $v
_{F}$ is also Henselian, whence $F _{\rm sep} \otimes _{F} F _{v
_{F}}$ is a field as well. As the exponent of Br$(F/E)$ divides
$[F\colon E]$, we have Br$(F) _{p} \neq \{0\}$, for each $p \in \Pi
(E)$. It is now easy to see from Lemma 3.2 that $F$ is a strictly
PQL-field.
\par
Step 2. Assume that finite extensions of $E$ are strictly
PQL-fields, and fix some $p \in \Pi (E)$. Clearly, if $M/E$ is a
finite Galois extension with $p \mid [M\colon E]$, and $M _{p}$ is
the fixed field of some Sylow $p$-subgroup of $\mathcal{G}(M/E)$,
then $p \in P(M _{p})$. Combining Lemmas 3.2 and 3.3 with (1.8), one
concludes that Br$(M _{p}) \cong {\rm Br}(E) _{p} \cong \ZZ (p
^{\infty })$. In particular, $p \in P(E)$, whence $E(p) \otimes _{E}
E _{v(p)}$ is a field, for some $v(p) \in M(E)$. Let now $F$ be a
finite Galois extension of $E$ in $E _{\rm sep}$ and $v(p) ^{\prime
}$ an absolute value of $F$ extending $v(p)$. In view of Lemma 1.3,
then Br$(F _{v(p)'}/E _{v(p)})$ is of exponent dividing $[F\colon
E]$, so it follows from Lemma 3.1 and the structure of Br$(E) _{p}$
that Br$(F _{v(p)'}) _{p} \neq \{0\}$ and the prolongation $v(p)
^{\prime }$ is unique. This means that $E _{\rm sep} \otimes _{E} E
_{v(p)}$ is a field (cf. \cite{CF}, Ch. II, Theorem 10.2), which
completes the proof of Corollary 3.4.
\end{proof}

\par
\medskip
\begin{rema} Lemma 3.3 indicates that if $E$ is an algebraic
strictly PQL-extension of a global field $E _{0}$, then finite
subgroups of Br$(E)$ are cyclic. This implies that the local
reciprocity law for arbitrary strictly PQL-fields (cf. \cite{Ch5},
Theorem 2) has the same form as in the case of local fields. In
other words, when $R/E$ is a finite abelian extension, $E ^{\ast
}/N(R/E) \cong \mathcal{G}(R/E)$.
\end{rema}

\par
\medskip
Note finally that if $E$ is an algebraic extension of a global field
$E _{0}$ satisfying the equivalent conditions of Corollary 3.4, then
$\mathcal{G} _{E}$ is prosolvable. This follows from (1.7) (iv) and
the solvability of the Galois groups of finite Galois extensions of
local fields (cf. \cite{FV}, Ch. IV, Sect. 1). When $E \in
I(\overline E _{0}/E _{0})$ is merely strictly PQL, $\mathcal{G}
_{E}$ is rarely prosolvable, as can be seen, for example, in the
formally real case, from Theorem 2.2 and the following result.
\par
\medskip
\begin{coro} Let $E$ be a formally real algebraic extension of $\QQ
$ admitting LCFT, and let $p \in P(E)$. Then there exists a finite
Galois extension of $E$ with a Galois group isomorphic to the
symmetric group {\rm Sym}$_{p}$.
\end{coro}

\par
\medskip
\begin{proof} It is clearly sufficient to consider only the special
case of $p \neq 2$. In this case, by (1.2), $E(p) \neq E$, i.e.
there exists a cyclic extension $L$ of $E$ of degree $p$. Fix $v(p)$
and $v(2)$ in $M(E)$ so that $E(p) \otimes _{E} E _{v(p)}$ and $E
(2) \otimes _{E} E _{v(2)}$ are fields. By Lemmas 3.2 and 3.3,
$v(p)$ and $v(2)$ are nonequivalent, so it follows from the (weak)
approximation theorem (cf. \cite{L2}, Ch. II, Theorem 1) and
Krasner's lemma that there exists a polynomial $f _{p} \in E[X]$ of
degree $p$ with the following properties: (i) the root field of $f
_{p} (X)$ over $E _{v(p)}$ is $E _{v(p)}$-isomorphic to $L \otimes
_{E} E _{v(p)}$; (ii) $f _{p} (X)$ has exactly $p - 2$ zeroes in $E
_{v(2)}$. This implies that the Galois group $G _{p}$ of $f _{p}
(X)$ over $E$ is isomorphic to a transitive subgroup of Sym$_{p}$
with a transposition, which ensures that $G _{p} \cong {\rm
Sym}_{p}$.
\end{proof}

\section{Minimal algebraic strictly PQL-extensions of global fields}

\par
\medskip
Let $E _{0}$ be a global field. Theorem 2.1 shows that a field $E
\in I(\overline E _{0}/E _{0})$ with $P(E) \neq \phi $ is strictly
PQL if and only if it possesses a system $V(E) = \{v(p) \in
M(E)\colon \ p \in P(E)\}$, such that $v(p _{i})$ is $p$-Henselian,
for each $p \in P(E)$, unless $p = 2 \in P(E)$ and $v(2)$ is
archimedean. Therefore, it is convenient to describe a number of
properties of the considered fields in the language of valuation
theory. Our objective in the present Section is to demonstrate this
by proving the existence of some families of algebraic strictly
PQL-extensions of $E _{0}$, in a form adequate to the needs of LCFT.
For the purpose, we need the following definition.
\par
\medskip
\begin{defn} Let $F$ be an algebraic extension of a global
field $E _{0}$ with Br$(F) \neq \{0\}$, $P$ a subset of $\{p \in \PP
\colon \ {\rm Br}(F) _{p} \neq \{0\}\}$, $P \neq \phi $, and $W =
\{w(p) \in M(F)\colon \ p \in P\}$ a system chosen so that Br$(F
_{w(p)}) _{p} \neq \{0\}$, for each $p \in P$. Assume that $\Pi =
\PP \setminus P$, $P(I) = P _{i}\colon \ i \in I$, is a partition of
$P$, i.e. a sequence of nonempty subsets of $P$, such that $\cup _{i
\in I} P _{i} = P$ and $P _{j} \cap P _{j'} = \phi $, for each pair
of distinct indices $j, j ^{\prime }$, $\chi = \{\chi _{i}\colon \ i
\in I\}$ is a system of group classes. The partition $P(I)$ is said
to be compatible with $W$, if the set $W _{i} = \{w(p _{i})\colon \
p _{i} \in P _{i}\}$ consists of equivalent absolute values; the
system $\chi $ is called admissible by the pair $(W, P(I))$, if the
following conditions hold:
\par
(c) $\chi _{i}$ consists of $P _{i} \cup \Pi _{i}$-groups, where
$\Pi _{i} \subseteq \Pi $, for each $i \in I$;
\par
(cc) $\chi _{i}$ is a saturated group formation including all finite
$p _{i}$-groups, for each $i \in I$ and $p _{i} \in P _{i}$, unless
$P _{i} = \{2\}$ and $w(2)$ is archimedean; if $P _{i} = \{2\}$ and
$w(2)$ is archimedean, then $\chi _{i}$ equals the class of groups
of orders $\le 2$.
\end{defn}

\par
\medskip
When $F$, $W$, $P$, $P(I)$ and $\chi $ are fixed as in Definition
4.1, we denote by $\Omega _{\chi } (F, P, W)$ the set of strictly
PQL-fields $K \in I(\overline E _{0}/F)$ with $P(K) = P$ and a
characteristic system $V(K) = \{v(p)\colon \ p \in P\}$, such that
$v(p)$ extends $w(p)$ whenever $p \in P$, the absolute values $v(p
_{i})$, $p _{i} \in P _{i}$, are equivalent, for each $i \in I$, and
$v(p _{i})$ is $\chi _{i}$-Henselian, in case $w(p _{i})$ is
nonarchimedean. Our main results concerning this set can be stated
as follows:
\par
\medskip
\begin{theo} Let $E _{0}$ be a global field and $F$ an
extension of $E _{0}$ in $\overline E _{0}$, such that {\rm Br}$(F)
\neq \{0\}$. Assume that $P$, $\Pi $ and $W = \{w(p)\colon \ p \in
P\}$ are defined as above, $P(I) = P _{i}\colon  i \in I$, is a
partition of $P$ compatible with $W$, and $\chi = \{\chi _{i}\colon
\ i \in I\}$ is a system of group classes admissible by $(W, P(I))$.
Then the set $\Omega _{\chi } (F, P, W)$ has the following
properties:
\par
{\rm (i)} $\Omega _{\chi } (F, P, W)$ is nonempty and satisfies the
conditions of Zorn's lemma with respect to the partial ordering
inverse to inclusion;
\par
{\rm (ii)} If $\chi _{i}$, $i \in I$, are abelian closed unless $P
_{i} = \{2\}$ and $w(2)$ is Archimedean, then every $E \in \Omega
_{\chi } (F, P, W)$ possesses a unique subfield $R(E)$, which is a
minimal element of $\Omega _{\chi } (F, P, W)$;
\par
{\rm (iii)} If $K \in \Omega _{\chi } (F, P, W)$ is minimal, $V(K) =
\{\kappa (p)\colon \ p \in P\}$ is a characteristic system of $K$,
and $F _{w(p)}$ is the closure of $F$ in $K _{\kappa (p)}$, then the
degrees of finite extensions of $F _{w(p)}$ in $K _{\kappa (p)}$ are
not divisible by $p$.
\end{theo}

\par
\medskip
{\it Proof.} We first show that $\Omega _{\chi } (F, P, W)$ contains
a minimal element $K$ with the properties required by Theorem 4.1
(iii). Denote by $R(p)$ the maximal $p$-extension of $R$ in
$\overline E _{0}$, provided that $R \in I(\bar F/F)$ and $p \in \PP
$. Fix an algebraic closure $\overline F _{w(p)}$ of $F _{w(p)}$ as
well as an embedding $\theta _{i}$ of $\overline F$ in $\overline F
_{w(p)}$ as an $F$-subalgebra, for each $i \in I$, and consider the
tower $\{K _{n}\colon \ n \in \NN \}$ of extensions of $F$ in
$\overline F$, defined inductively as follows:
\par
\medskip
(4.1) $K _{1} = F$, and for each $n \in \NN $, $K _{n+1}$ is the
compositum of the fields $K _{n,i} = \{f _{i} \in K _{n} (\chi
_{i})\colon  \theta _{i} (f _{i}) \in \theta _{i} (K _{n})F
_{w(p)}\}$, $i \in I$, and $K _{n} (\pi )$, $\pi \in \Pi $.
\par
\medskip
Denote by $K$ the union $\cup _{n=1} ^{\infty } K _{n}$, and put
$\kappa (p _{i}) (\alpha ) = \bar w(p _{i}) (\theta _{i} (\alpha
))$, for each $\alpha \in K$, in case $p _{i} \in P _{i}$ and $i \in
I$, where $\bar w(p _{i})$ is the unique absolute value of
$\overline F _{w(p _{i})}$ extending the continuous prolongation of
$w(p _{i})$ on $F _{w(p _{i})}$. Also, let $\kappa _{\nu } (p)$ be
the absolute value of $K _{\nu }$ induced by $\kappa (p)$, for each
pair $(p, \nu ) \in (P \times \NN )$. We show that $K$ is a minimal
element of $\Omega _{\chi } (F, P, W)$ with a characteristic system
$V(K) = \{\kappa (p)\colon  p \in P\}$. It follows from (1.1) (ii)
that every finite Galois extension $L$ of $K$ in $\overline F$ is
isomorphic over $K$ to $L _{n} \otimes _{K _{n}} K$, for some index
$n$ depending on $L$, and a suitably chosen finite Galois extension
$L _{n}$ of $K _{n}$ in $L$. Since $\mathcal{G}(L _{n}/K _{n}) \cong
\mathcal{G}(L/K)$, this ensures that $L _{n} \cap K = L _{n} \cap K
_{n+1} = K _{n}$. It is therefore clear that if $p \in P _{i'}$ and
$\mathcal{G}(L _{n}/K _{n}) \in \chi _{i'}$, for some $i ^{\prime }
\in I$, then $\kappa _{n} (p)$ and $\kappa (p)$ have unique
prolongations on $L _{n}$ and $L$, respectively. This means that if
$w(p)$ is nonarchimedean, then $\kappa (p)$ is $\chi
_{i'}$-Henselian. A similar argument also shows that $K(p) = K$, in
case $p \in \Pi $. Now fix an index $i \in I$, and a prime $p _{i}
\in P _{i}$, and for each $R \in I(\overline E _{0}/F)$, let $w _{R}
(p _{i})$ be the absolute value of $\theta _{i} (R).F _{w(p _{i})}$
induced by $\bar w(p _{i})$, and $\widetilde R _{p _{i}}$ the
completion of $\theta _{i} (R).F _{w(p _{i})}$ with respect to $w
_{R} (p _{i})$. Then $w _{R} (p _{i})$ is Henselian, which means
that $\theta _{i} (R).F _{w(p _{i})}$ is separably closed in
$\widetilde R _{p _{i}}$. In view of (4.1) and Definition 4.1, this
result, applied to the case of $R = K$, implies that $p _{i}$ does
not divide the degrees of finite extensions of $F _{w(p _{i})}$ in
$\widetilde K _{p _{i}}$. Hence, by (1.3) (ii) and the behaviour of
Schur indices of central simple algebras under finite extensions of
their centres (cf. [P, Sect. 13.4]), we have Br$(\widetilde K _{p
_{i}}) _{p _{i}} \neq \{0\}$. Observing also that there exists an
$F$-isomorphism $K _{\kappa (p _{i})} \cong \widetilde K _{p _{i}}$
acting on $K$ as $\theta _{i}$, one concludes that Br$(K _{\kappa (p
_{i})}) _{p _{i}} \neq \{0\}$ and $p _{i}$ does not divide the
degrees of finite extensions of $F _{w(p)}$ in $K _{\kappa (p
_{i})}$. This, combined with (1.2) and Proposition 1.2, yields
consecutively Br$(K) _{p _{i}} \neq \{0\}$ and $K(p _{i}) \neq K$;
one also sees that if $P _{i} = \{2\}$ and $w(2)$ is archimedean,
then $K$ is formally real and $K _{\kappa (2)} \cong \RR $. The
assumptions of Theorem 4.1 and the obtained results indicate that
$K(p) \otimes _{K} K _{\kappa (p)}$ is a field, for each $p \in P$,
which completes the proof of the assertion that $K \in \Omega (F, P,
W)$ and $V$ is a characteristic system of $K$. Since $\kappa (p
_{i})$, $p _{i} \in P _{i}$, are Henselian whenever $P _{i}$
contains at least $2$ elements, these results enable one to deduce
from Lemma 3.1 (ii) or Lemma 1.1 that $\kappa (p _{i})$, $p _{i} \in
P _{i}$, are equivalent, whence $K \in \Omega _{\chi }(F, P, W)$. It
remains to be seen that $K$ is a minimal element of $\Omega _{\chi }
(F, P, W)$. Let $M$ be a proper extension of $F$ in $K$. Then there
exists an index $m \in \NN $, such that $K _{m} \subseteq M$ and $K
_{m+1} \not\subseteq M$. This means that $K _{m} (\pi )
\not\subseteq M$, for some $\pi \in \Pi $, or $K _{m,j}
\not\subseteq M$, for some $j \in I$. In the former case, this means
that $MK _{m} (\pi )$ contains as a subfield a cyclic extension of
$M$ of degree $\pi $, whence $\pi \in P(M)$. In the latter one, it
turns out that $M$ admits a proper finite extension $M ^{\prime }$
in $MK _{m,j}$. Let $M _{1}$ be the normal closure of $M ^{\prime }$
in $\overline E _{0}$ over $M$. Then $M _{1}$ is a subfield of $MK
_{m} (\chi _{i})$, so it follows from Galois theory that
$\mathcal{G}(M _{1}/M) \in \chi _{i}$. On the other hand, (4.1), the
inequality $M ^{\prime } \neq M$ and the inclusions $M \subseteq M
^{\prime } \subseteq MK _{m,j}$ indicate that the absolute value of
$M$ induced by $\kappa (p _{j})$, for a given $p _{j} \in P _{j}$,
has at least $2$ different prolongations on $M ^{\prime }$ and on $M
_{1}$. Applying now Lemma 3.1, one sees that $M \not\in \Omega
_{\chi } (F, P, W)$, i.e. $K$ is minimal in $\Omega _{\chi } (F, P,
W)$, as claimed.
\par
In the rest of the proof of Theorem 4.1 we identify the completions
$E _{v(p _{i})}$, $p _{i} \in P _{i}$, for each $i \in I$. Our
considerations rely on the following lemma.
\par
\medskip
\begin{lemm} Let $E _{0}$ be a global field, $E$ and $T$ extensions of
$E _{0}$ in $\overline E _{0}$, such that $T \subseteq E$, and let
$v \in M(E)$ and $u$ be the absolute value of $T$ induced by $v$.
Assume that $E$ and $T$ admit local $p$-class field theory, for some
$p \in P(E)$, and $E(p) \otimes _{E} E _{v}$ is a field. Then $T (p)
\otimes _{T} T _{u}$ is a field and $p \in P(T)$.
\end{lemm}

\par
\medskip
\begin{proof} Since $T _{u}$ is isomorphic to the closure of $T$ in $E
_{v}$, this can be deduced from (1.8) and Lemmas 3.1 (ii) and 3.2.
\end{proof}
\par
\medskip
We turn to the proof of Theorem 4.1 (i). Fix an element $E$ of
$\Omega _{\chi } (F, P, W)$, a characteristic system $V(E) =
\{v(p)\colon \ p \in P\}$ of $E$, and a nonempty linearly ordered
subset $\Lambda $ of $I(E/F) \cap \Omega _{\chi } (F, P, W)$. Denote
by $L$ the intersection of the fields from $\Lambda $, and by $v
_{L} (p)$ the absolute value of $L$ induced by $v(p)$, for each $p
\in P$. We show that $L \in \Omega _{\chi } (F, P, W)$ and $V(L) =
\{v _{L} (p)\colon \ p \in P\}$ is a characteristic system for $L$.
Note first that $P(L) = P$. Indeed, (1.2), (1.8) and Theorem 2.1
indicate that Br$(L) _{p} \neq \{0\}$ and $L(p) \neq L$, for every
$p \in P$. Assuming that $\tilde p \in \PP \setminus P$, one obtains
from Galois theory that $L(\tilde p).R/R$ is a $\tilde p$-extension,
which yields $L(\tilde p).R = R$, for each $R \in \Lambda $, and so
proves that $L (\tilde p) = L$. Suppose now that $i \in I$, $p _{i}
\in P _{i}$, $R \in \Lambda $, $v _{R} (p _{i})$ is the absolute
value of $R$ induced by $v(p _{i})$, and $\Sigma _{i}$ is an
algebraic closure of $E _{v(p _{i})}$. Identifying $\overline E
_{0}$ with its $L$-isomorphic copy in $\Sigma _{i}$, and $L _{v _{L}
(p)}$ and $R _{v _{R} (p _{i})}$ with the closures in $E _{v(p
_{i})}$ of $L$ and $R$, respectively, one obtains from Galois
theory, general properties of tensor products (cf. [P, Sect. 9.2,
Proposition c]), the choice of $v _{R} (p _{i})$, and Lemma 4.2 that
$R(\chi _{i}) \cap R _{v _{R} (p _{i})} = R$. Since $L (\chi _{i}).R
\subseteq R (\chi _{i})$, this implies that $L (\chi _{i}) \cap L
_{v _{L} (p _{i})} \subseteq R$. As $R$ is an arbitrary element of
$\Lambda $, the obtained result proves that $L (\chi _{i}) \cap L
_{v _{L} (p _{i})} = L$. This means that $L (\chi _{i}) \otimes _{L}
L _{v _{L} (p _{i})}$ is a field, so Theorem 4.1 (i) is proved.
\par
Our objective now is to prove Theorem 4.1 (iii). Let $L$ be a
minimal element of $\Omega _{\chi }(F, P, W)$ and $V(L) = \{\lambda
(p)\colon \ p \in P\}$. It is clearly sufficient to show that $L$
can be viewed as an extension of $F$ defined in accordance with
(4.1). Identifying, for each $i \in I$ and $p _{i} \in P _{i}$, $L$
with its canonically isomorphic copy in $L _{\lambda (p _{i})}$, and
$F _{w(p _{i})}$ with the closure of $F$ in $L _{\lambda (p _{i})}$,
fix an embedding $j _{i}$ of $\overline L$ in $\Sigma
_{i}$ as an $E$-subalgebra. Also, let $K$ be the union of the
extensions $K _{n}\colon \ n \in \NN $, of $F$ in $\overline E$
associated with the maps $j _{i}\colon  i \in I$, as in (4.1).
Proceeding by induction on $n$, and arguing as in the proof of
Theorem 4.1 (i), one obtains that $K _{n} \subseteq E$, for every $n
\in \NN $, i.e. $K \subseteq E$. Since $K$ is a minimal element of
$\Omega _{\chi }(F, P, W)$, this yields $K = E$, and so completes
the proof of Theorem 4.1 (iii).
\par
It remains for us to prove Theorem 4.1 (ii). In view of Theorem 4.1
(i), it suffices to show that if $E \in \Omega _{\chi }(F, P, W)$,
then $I _{\chi }(E/F)$ contains a unique minimal element. Assuming
the opposite, take a field $E \in \Omega _{\chi }(F, P, W)$ so that
$I _{\chi }(E/F)$ contains two different elements, say, $K$ and $L$.
As shown in the proof of Theorem 4.1 (iii), $K = \cup _{n=1}
^{\infty } K _{n}$ and $L = \cup _{n=1} ^{\infty } L _{n}$, the
unions being defined as in (4.1). As $K \neq L$, this means that $K
_{n} = L _{n}$ and $K _{n+1} \neq L _{n+1}$, for some $n \in \NN $.
Specifically, $K _{n,i} \neq L _{n,i}$, for some $i \in I$. In view
of (1.6), this means that $R _{i}(\chi _{i}) = R _{i}$, where $R
_{i} = K _{n,i}L _{n,i}$. Hence, by Definition 3 and the assumptions
of Theorem 4.1 (ii), $p _{i} \notin P(R _{i})$, for any $p _{i} \in
P _{i}$. It is therefore clear from (1.8) that Br$(R _{i} ^{\prime
}) _{p _{i}} = \{0\}$, for every $R _{i} ^{\prime } \in I(\overline
E _{0}/R _{i})$ and each $p _{i} \in P _{i}$. Since $E \in
I(\overline E _{0}/R _{i})$ and Br$(E) _{p _{i}} = \{0\}$ when $p
_{i}$ runs across $P _{i}$, this is a contradiction proving our
assertion about $I _{\chi }(E/F)$, the concluding step towards the
proof of Theorem 4.1.
\par
\medskip
{\it Proof of Theorem 2.2.} Statements (i), (iii) and (iv) are
obtained by applying Theorem 4.1 to the special case in which $I =
P$ and, for each $p \in P$, $\chi _{p}$ is the formation of finite
$p$-groups unless $2 \in P$ and $w(2)$ is archimedean. Since the
class of solvable groups is closed under taking subgroups, quotient
groups and group extensions, this enables one to deduce Theorem 2.2
(ii) from Galois theory and (4.1).
\par
\medskip
\begin{coro} Let $E _{0}$ be a global field, $F$ an algebraic
extension of $E _{0}$, such that Br$(F) \neq \{0\}$, and $P$ a
nonempty subset of $\PP $, for which {\rm Br}$(F) _{p} \neq \{0\}$,
$p \in P$. Then the set $\Omega _{P} (\overline E _{0}/F)$ of
strictly {\rm PQL}-fields $\Sigma \in I(\overline E _{0}/F)$ with
$P(\Sigma ) = P$ is nonempty. Moreover, every $K \in \Omega _{P}
(\overline E _{0}/F)$ possesses a unique subfield $R(K)$ that is a
minimal element of $\Omega _{P} (\overline E _{0}/F)$.
\end{coro}

\par
\medskip
\begin{proof} It is easily obtained from (1.8) and Lemma 4.2 that
$\Omega _{P} (\overline E _{0}/F)$ equals the union of the sets
$\Omega (F, P, W)$, taken over all $W$ admissible by Definition 1.
Observing also that $I(E/F) \cap \Omega _{P} (\overline E _{0}/F)
\subseteq \Omega (F, P, W)$, for each $E \in \Omega (F, P, W)$, one
proves our assertions by applying Theorem 4.1 to the set $\Omega
_{\chi } (F, P, W)$ considered in the proof of Theorem 2.2, for an
arbitrary system $W = \{w(p) \in M(F)\colon \ p \in P\}$, such that
Br$(F _{w(p)}) _{p} \neq \{0\}$, $p \in P$.
\end{proof}

\par
\medskip
\begin{coro} Under the hypotheses of Theorem 4.1,
suppose that $F = E _{0}$, $2 \in P$, $w(2)$ is Archimedean, $K$ is
a minimal element of $\Omega _{\chi } (F, P, W)$, and $n$ is an
integer $\ge 5$. Then $I(\overline E _{0}/K)$ contains Galois
extensions $M _{k}, L _{k}$, $k \in \NN $, of $K$, such that
$\mathcal{G}(M _{k}/K) \cong {\rm Sym}_{n}$ and $\mathcal{G}(L
_{k}/K) \cong {\rm Alt}_{n}$, for every $k \in \NN $.
\end{coro}

\par
\medskip
\begin{proof} Let $d _{1}$ and $d _{2}$ be squarefree integers, such
that $d _{1} < 0 < d _{2}$ and $d _{2} \in E _{0} ^{\ast 2}$ (the
existence of $d _{2}$ is implied Lemma 1.1). It is known (cf.
\cite{G}) that there exists a set $\{A _{k}, B _{k}\colon \ k \in
\NN \}$ of Galois extensions of $\QQ $ in $\overline F$ with
$\mathcal{G}(A _{k}/\QQ ) \cong \mathcal{G}(B _{k}/\QQ ) \cong {\rm
Sym}_{n}$, $\sqrt{d _{1}} \in A _{k}$ and $\sqrt{d _{2}} \in B
_{k}$, for every $k \in \NN $. Note also that the Galois closure $K
^{\prime }$ of $K$ in $\overline F$ over $F$ has a prosolvable
Galois group. Since, by the proof of Theorem 4.1, $K$ is obtained
from $F$ in accordance with (4.1), and the class of solvable groups
is closed under taking subgroups, quotient groups and group
extensions, the prosolvability of $\mathcal{G}(K ^{\prime }/F)$ can
be deduced from Galois theory and the solvability of finite groups
of odd order [FT]. Hence, by Galois theory and the simplicity of
Alt$_{n}$, $n \ge 5$, when $k$ is sufficiently large, the fields $A
_{k}K$ and $B _{k}K$ are Galois extensions of $K$ with
$\mathcal{G}(A _{k}K/K) \cong {\rm Sym}_{n}$ and $\mathcal{G}(B
_{k}K/K) \cong {\rm Alt}_{n}$.
\end{proof}

\begin{rema}
Corollary 4.4 retains validity, if $F = E _{0}$, $2 \in P$ and, in
terms of Definition 3, $P _{j} = \{2\}$ and $\Pi _{j}$ contains at
most one element, for some $j \in I$; in particular, this holds when
$\PP \setminus P$ contains at most one element and $w(2)$ is
nonequivalent to $w(p)$, for any $p \in P \setminus \{2\}$. Indeed,
in this case, $\mathcal{G}(K ^{\prime }/F)$ is prosolvable, $K
^{\prime }$ being again the Galois closure of $K$ in $\overline E
_{0}$ over $F$ (as follows from the solvability of finite groups of
odd or biprimary order). Omitting the details, note that the
extensions $A _{k}$ and $B _{k}$ of $F$ can be constructed as root
fields of polynomials in $F[X]$ of degree $n$ with a suitably chosen
local behaviour (ensured by applying the approximation theorem).
\end{rema}

\medskip
\begin{prop} In the setting of Theorem 4.1, let $F$ be a global
field. Then the minimal elements of $\Omega _{\chi } (F, P, W)$ form
a single $F$-isomorphism class if and only if some of the following
assertions holds true:
\par
{\rm (i)} $P = P _{1}$;
\par
{\rm (ii)} $P = \PP $ and for each $i \in I$, $p _{i} \in P _{i}$,
the system $W _{i} = \{w(q _{i})\colon \ q _{i} \in P _{i}\}$
contains all $w(p) \in W$ that are equivalent to $w(p _{i})$.
\end{prop}

\begin{proof} We first show that if none of conditions (i) and (ii)
holds, then $\Omega _{\chi } (F, P, W)$ contains a pair of
nonisomorphic fields (over $F$). Our argument relies on two
observations the latter of which is a special case of Lemma 1.1:
\par
\medskip
(4.2) (i) If $E _{1}$ and $E _{2}$ are fields lying in $\Omega
_{\chi } (F, P, W)$, with characteristic systems $V(E _{u}) = \{v
_{u} (p)\colon \ p \in P\}$, $u = 1, 2$, and if there exists an
$F$-isomorphism $\varphi $: $E _{1} \to E _{2}$, then the absolute
values $v _{1} (p)$ and $v _{2} (p) \circ \varphi $ of $E _{1}$ are
equivalent, for each $p \in P$;
\par
(ii) If $\pi \in \PP $, $p \in P$ and $S$ is a finite subset of $W$,
then there exists a cyclic extension $\Phi _{\pi }$ of $F$ in $F
_{{\rm sep}}$, such that $[\Phi _{\pi }\colon F] = \pi $ and each $s
\in S$ has $\pi $ distinct prolongations on $\Phi _{\pi }$.
\par
\medskip
Suppose first that $\pi \not\in P$ and there exist $p _{1}, p _{2}
\in P$, for which $w(p _{1})$ is not equivalent to $w(p _{2})$.
Assume also that $S = \{w(p _{1}), w(p _{2})\}$, fix a generator
$\psi _{\pi }$ of $\mathcal{G}(\Phi _{\pi }/F)$ and prolongations
$w(p _{1}) ^{\prime }$ and $w(p _{2} ^{\prime })$ on $\Phi _{\pi }$
of $w(p _{1})$ and $w(p _{2})$, respectively. By Lemma 1.3, the set
$S _{u}$ of absolute values of $\Phi _{\pi }$ extending $w(p _{u})$
equals $\{w(p _{u}) ^{\prime } \circ \psi _{\pi } ^{j}\colon  j = 0,
1,..., \pi _{\pi } - 1\}$, for each index $u$. It is therefore clear
that $\Phi _{\pi } \in I(E/F)$, for any $E \in \Omega _{\chi } (F,
P, W)$. Put $W _{u} = \{w _{u} (p)\colon  p \in P\}$, where $w _{u}
(p)$ is a prolongation of $w(p)$ on $\Phi _{\pi }$, for each $p \in
P$, $w _{1} (p _{u}) = w(p _{u}) ^{\prime }$ and $w _{2} (p _{u}) =
w(p _{u}) ^{\prime } \circ \psi _{\pi } ^{u-1}\colon \ u = 1, 2$.
Applying (4.1) and Lemma 1.3, one obtains that $\Omega _{\chi } (F,
P, w)$ contains minimal elements $E _{1}$ and $E _{2}$ possessing
characteristic systems $V(E _{1})$ and $V(E _{2})$, respectively,
such that $V(E _{u}) = \{v _{u} (p)\colon \ p \in P\}$ and $v
_{u}(p)$ extends $w _{u} (p)$ whenever $p \in P$ and $u = 1, 2$.
Hence, by (4.2) (i) and the choice of $w _{u} (p _{1})$ and $w _{u}
(p _{2})$, the fields $E _{1}$ and $E _{2}$ are not isomorphic over
$F$.
\par
Suppose now that there exist different indices $j$, $m \in I$, for
which $\{w(p ^{\prime })\colon $ $p ^{\prime } \in P _{j} \cup P
_{m}\}$ consists of equivalent absolute values. Assume also that
$\pi \in P _{j}$ and $S = \{w(\pi )\}$, and fix some $p \in P _{m}$.
It is easily verified that the field $\Phi _{\pi }$ considered in
(4.2) (ii) is included in all fields from $\Omega _{\chi } (F, P,
W)$. Arguing as in the case of $\pi \notin P$, one obtains that
$\Omega _{\chi } (F, P, W)$ contains minimal elements $E _{1}$ and
$E _{2}$, such that the system $\{v _{1} (p ^{\prime }) \in V(E
_{1})\colon \ p ^{\prime } \in P _{1} \cup P _{m}\}$ consists of
equivalent absolute values, whereas $\{v _{2} (p ^{\prime }) \in V(E
_{2})\colon \ p ^{\prime } \in P _{i} \cup P _{m}\}$ does not
possess this property. It is therefore clear from (4.2) (i) that $E
_{1}$ and $E _{2}$ are not $F$-isomorphic, which completes the proof
of the left-to-right implication of Proposition 4.6.
\par
Conversely, let some of conditions (i), (ii) be in force, take
minimal elements $K$ and $L$ of $\Omega _{\chi } (F, P, W)$, as well
as characteristic systems $V(K) = \{\kappa (p)\colon \ p \in P\}$
and $V(L) = \{\lambda (p)\colon \ p \in P\}$, and for each $n \in
\NN $, denote by $K _{n}$ and $L _{n}$ the fields in $I(K/F)$ and
$I(L/F)$, respectively, defined in accordance with (4.1). Also, let
$\kappa _{n} (p)$ and $\lambda _{n} (p)$ be the absolute values of
$K _{n}$ and $L _{n}$, induced by $\kappa (p)$ and $\lambda (p)$,
respectively, for any index $n$. Proceeding by induction on $n$ (and
taking the inductive step by a repeated application of Lemma 1.3),
one proves the existence of a sequence $\psi = \{\psi _{n}\colon \ n
\in \NN \}$ of $F$-isomorphisms $\psi _{n}\colon \ K _{n} \to L
_{n}$, such that $\psi _{n+1}$ extends $\psi _{n}$ and $\kappa _{n}
(p)$ equals the composition $\lambda _{n} (p) \circ \psi _{n}$, for
each pair $(n, p) \in \NN \times P$. The sequence $\psi $ gives rise
to an $F$-isomorphism $K \cong L$, so Proposition 4.6 is proved.
\end{proof}

\par
\medskip
Proposition 4.6 shows that if $F$ is a global field, then the
minimal elements of $\Omega (F, P, W)$ form an $F$-isomorphism class
if and only if $P$ contains only one element or $P = \PP $ and $W$
consists of pairwise nonequivalent absolute values. The concluding
result of this Section determines the structure of the Brauer groups
of algebraic strictly PQL-extensions of global fields, and gives a
classification of abelian torsion groups realizable as such Brauer
groups.
\par
\medskip
\begin{prop} Let $E _{0}$ be a global field, $T$ be a nontrivial
divisible abelian torsion group, $P(T)$ the set of those $p \in \PP
$ for which $T _{p} \neq \{0\}$, and $T _{0}$ the direct sum $\oplus
_{p \in P _{0} (T)} \ZZ (p ^{\infty })$, where $P _{0} (T) = P(T)
\setminus \{2\}$. Then $T \cong {\rm Br}(E(T))$, for some strictly
{\rm PQL}-field $E(T) \in I(\overline E _{0}/E _{0})$, if and only
if $T$ is isomorphic to one of the following groups:
\par
{\rm (i)} the direct sum of $T _{0}$ by a group of order $2$; when
this occurs, $E(T)$ is formally real;
\par
{\rm (ii)} $\oplus _{p \in P(T)} \ZZ (p ^{\infty })$; in this case,
$E(T)$ is nonreal and can be chosen so that its finite extensions
are strictly {\rm PQL}.
\end{prop}
\par
\medskip
\begin{proof} Let $E$ be an extension of $E _{0}$ in $\overline E
_{0}$. It is clear from (1.2) that Br$(E) _{q} = \{0\}$, for every
$q \in \PP \setminus P(E)$. This, combined with Lemma 3.3 and
Theorems 2.1 and 2.2, proves all assertions of Proposition 4.7
except the one that if $T$ is of type (ii), then $E(T)$ can be
chosen so that its finite extensions are strictly PQL. At the same
time, Theorem 4.1, applied to a saturated group formation $\chi $
(including the class of finite $p$-groups, for every $p \in P(T)$),
the partition $P = P _{1} = P(T)$, and a nonarchimedean absolute
value $w _{1}$ of $E _{0}$, indicates that there exists a $\chi
$-Henselian algebraic strictly PQL-extension $E$ of $E _{0}$ with
$P(E) = P(T)$. Hence, by Lemma 3.3, Br$(E) \cong T$ if and only if
$T$ is of type (ii). By Corollary 3.4, when $\chi $ is the formation
of all finite groups, finite extensions of $E$ are strictly PQL, so
the obtained result completes the proof of Proposition 4.7.
\end{proof}

\section{Hasse norm principle}

\par
\medskip
In this Section we prove the validity of the Hasse norm principle
for arbitrary finite Galois extensions of the fields pointed out at
the end of the Introduction. By Theorem 2.1, this applies to
algebraic strictly PQL-extensions of global fields, which is used
for proving Theorem 2.3 and Proposition 2.4.
\par
\medskip
\begin{theo} Assume that $E _{0}$ is a global field and $E$ is an
extension of $E _{0}$ in $\overline E _{0}$, such that the set $S
_{p} (E) = \{v \in M(E)\colon \ {\rm Br}(E _{v}) _{p} \neq \{0\}\}$
is finite, for each $p \in \PP $. Let $M/E$ be a finite Galois
extension. Then:
\par
{\rm (i)} $N(M/E)$ contains every $\lambda  \in E ^{\ast }$
satisfying the inequalities $v(\lambda - 1) < \varepsilon \colon \ v
\in \cup _{p \in P(M/E)} S _{p} (E)$, for a some real number
$\varepsilon
> 0$ depending on $M/E$;
\par
{\rm (ii)} $N(M/E) = N _{\rm loc} (M/E)$; an element $c \in E ^{\ast }$
lies in $N(M/E)$ if and only if $c \in N(M _{v'}/E _{v})$ whenever
$v \in \cup _{p \in P(M/E)} S _{p} (E)$ and $v ^{\prime }$ is an
absolute value of $M$ extending $v$.
\end{theo}

\par
\medskip
{\it Proof.} We first show that $N(M/E) \subseteq N _{\rm loc}
(M/E)$, i.e. of $N(M _{v'}/E _{v})$, for $v \in E$, and an arbitrary
prolongation $v ^{\prime }$ of $v$ on $M$. Identifying $E _{v}$ with
the topological closure of $E$ in $M _{v'}$, one obtains without
difficulty that $M _{v'}/E _{v}$ is a Galois extension with
$\mathcal{G}(M _{v'}/E _{v}) \cong \mathcal{G}(M/(M \cap E _{v}))$.
Observing also that $N _{E _{v}} ^{M ^{v'}} (\mu ) = N _{F} ^{M}
(\mu )$, for every $\mu \in M ^{\ast }$, where $F = M \cap E _{v}$,
one deduces the required inclusion from the following lemma proved
in \cite{Ch2}, II.

\par
\medskip
\begin{lemm} Assume that $E$, $M$ and $F$ are fields, such that
$M/E$ is a finite Galois extension and $F \in I(M/E)$. For each $p
\in P(M/E)$, let $L _{p}$ be the fixed field of some Sylow
$p$-subgroup $G _{p}$ of $\mathcal{G}(M/E)$. Then $N(M/E) \subseteq
N(M/F)$ and $N(M/E) = \cap _{p \in P(M/E)} N(M/L _{p})$.
\end{lemm}

\par
\medskip
Note that the class of fields satisfying the conditions of Theorem
5.1 is closed under the formation of finite extensions. This follows
from (1.8) and the fact that absolute values of fields have finitely
many prolongations on their finite extensions. Hence, by Lemma 5.2,
it suffices to establish the rest of Theorem 5.1 in the special case
where $M \subseteq E(p)$.
\par
\medskip
\begin{lemm} Under the hypotheses of Theorem 5.1, let
$M/E$ be a finite $p$-extension, for some $p \in \PP $. Then:
\par
{\rm (i)} An element $c \in E ^{\ast }$ lies in $N(M/E)$, provided
that $c \in N(M _{v'}/E _{v})$ when $v$ runs across $S _{p} (E)$ and
$v ^{\prime }$ is a prolongation of $v$ on $M$;
\par
{\rm (ii)} There exists a real number $\varepsilon > 0$, such that
$N(M/E)$ contains every $\lambda \in E ^{\ast }$ satisfying
$v(\lambda - 1) < \varepsilon \colon \ v \in S _{p} (E)$.
\end{lemm}

\par
\medskip
\begin{proof} For each $v \in S _{p} (E)$ and any absolute value $v
^{\prime }$ of $M$ extending $v$ , denote by $\bar v$ and $\bar v
^{\prime }$ their continuous prolongations on $E _{v}$ and $M _{v'}$
respectively. Also, let $\eta _{v'}$ be a primitive element of $M
_{v'}/E _{v}$ taken so that $\bar v ^{\prime } (\eta _{v'} - 1) <
1$, and let $h _{v'} (X)$ be the minimal polynomial of $\eta _{v'}$
over $E _{v'}$. It follows from Krasner's lemma that $M _{v'}/E
_{v}$ contains as a primitive element a root of the polynomial $h
_{v'} (X) + (\lambda _{v} - 1).a _{v'}$, provided that $a _{v'}$ is
the free term of $h _{v'} (X)$ and $\lambda _{v}$ is any element of
$E _{v} ^{\ast }$ for which $\bar v (\lambda _{v} - 1)$ is
sufficiently small. As $S _{p} (E)$ is finite, this enables one to
deduce Lemma 5.3 (ii) from  Lemma 5.3 (i) and \cite{Ch2}, I, Lemma
4.2 (ii).
\par
We prove Lemma 5.3 (i). Suppose that $[M\colon E] = p ^{n}$, for
some $n \in \NN $, fix a field $F \in I(M/E)$ so that $[F\colon E] =
p$, and denote by $\varphi $ some generator of $\mathcal{G}(F/E)$.
Assume that $c \in N _{\rm loc} (M/E)$. Considering the cyclic
$E$-algebra $(F/E, \varphi , c)$, one obtains from Proposition 1.2
and \cite{P}, Sect. 15.1, Proposition b, that $c \in N(F/E)$, which
proves Lemma 5.3 in case $n = 1$.
\par
Let now $n \ge 2$ and $\alpha  \in F ^{\ast }$ be an element of norm
$N _{E}^{F}(\alpha ) = c$. We show that there exists $\beta \in F
^{\ast }$, such that $\alpha (\varphi (\beta )\beta ^{-1}) \in
N(M/F)$. Denote by $S(F)$ the set of absolute values of $F$
extending elements of $M(E)$ and by $\tilde S _{p} (F)$ the subset
of $S(F)$ consisting of prolongations of elements of $S _{p} (E)$.
It is clear from (1.8) that $S _{p} (F) \subseteq \tilde S _{p} (F)$
and equality holds in case $E$ is a nonreal field. Proceeding by
induction on $n$, one concludes that it suffices to establish the
existence of $\beta $, under the hypothesis that $N(M/F)$ has the
properties required by Lemma 5.3. Fix an element $\omega $ of $S
_{p} (F)$ as well as a prolongation $v ^{\prime }$ of $\omega $ on
$M$, and denote by $v$ the absolute value of $E$ induced by $\omega
$. Suppose first that $\omega $ is the unique prolongation of $v$ on
$F$. Then $F _{\omega }$ is $E _{v}$-isomorphic to $F \otimes _{E} E
_{v}$. This means that $F _{\omega }/E _{v}$ is a cyclic extension
of degree $p$, or equivalently, that $\varphi $ extends uniquely to
an $E _{v}$-automorphism $\varphi _{\omega }$ of $F _{\omega }$. As
$c \in N(R _{v'}/E _{v})$, $N _{E} ^{F} (\alpha ) = c$ and $E$ is
dense in $E _{v}$, the obtained result, combined with Hilbert's
Theorem 90 and the inductive hypothesis, proves that $\alpha \varphi
(\beta _{v'})\beta _{v'} ^{-1} \in N(R _{v'}/F _{\omega })$, for
some $\beta _{v'} \in F ^{\ast }$.
\par
Suppose finally that $v$ is not uniquely extendable on $F$, and fix
a real number $\varepsilon > 0$. Then, by Lemma 1.3, the
compositions $\omega \circ \varphi ^{i} := \omega _{i}$, $i = 0,...,
p - 1$, are the prolongations of $v$ on $F$ (and are pairwise
nonequivalent). Note also that $E$ is a dense subfield of $F
_{\omega _{0}},..., F _{\omega _{p-1}}$, so one can find elements $c
_{0},..., c _{p-1}$ of $E ^{\ast }$ so that $\omega _{i} (\alpha - c
_{i}) < \varepsilon $, for each index $i$. Clearly, $E ^{\ast }$
contains elements $\mu _{0},..., \mu _{p-1}$ satisfying the
equalities $\mu _{j}\mu _{j+1} ^{-1} = c _{j}\colon \ j = 0,..., p -
2$, and by the approximation theorem, there exists $\beta _{\omega}
\in F ^{\ast }$, such that $\omega _{j} (\beta _{\omega } - \mu
_{j}) < \varepsilon $, for every $j \in \{0,..., p - 1\}$.
Therefore, the product $\alpha _{\omega } := \alpha \varphi ^{(p-1)}
(\beta _{\omega }).\beta _{\omega } ^{-1}$ satisfies the
inequalities $\omega _{j} (\alpha _{\omega } - 1) < \varepsilon $,
$j \in \{0,..., p - 2\}$, and $\omega _{p-1} (\alpha _{\omega } - c)
< \varepsilon $. When $\varepsilon $ is sufficiently small, this
ensures that $\alpha _{\omega } \in N(M _{v _{i}}/F _{\omega
_{i}})$, provided that $i \in \{0,1,..., p - 1\}$ and $v _{i}$ is an
absolute value of $M$ extending $\omega _{i}$ (cf. \cite{L2}, Ch.
II, Proposition 2). Summing up the obtained results and using again
the approximation theorem, one concludes that there exists $\beta
\in F ^{\ast }$ such that $\alpha \varphi (\beta )\beta ^{-1} \in
N(M _{\omega '}/F _{\omega })$, for each $\omega \in S _{p} (F)$,
and any prolongation $\omega ^{\prime }$ of $\omega $ on $M$. Now
the inductive hypothesis leads to the conclusion that $\alpha
\varphi (\beta )\beta ^{-1} \in N(M/F)$ and $c \in N(M/E)$, which
completes the proof of Lemma 5.3 and Theorem 5.1.
\end{proof}

\par
\medskip
\begin{coro} Assume that $E _{0}$, $E$ and $S _{p} (E)$, $p \in \PP
$, satisfy the conditions of Theorem 5.1, $R$ is a finite extension
of $E$ in $E _{\rm sep}$, $M$ is the normal closure of $R$ in $E
_{\rm sep}$ over $E$, $f$ is a norm form of $R/E$, and $c$ is an
element of $E ^{\ast }$. Then $c \in N(R/E)$ if and only if $c$ is
presentable by $f$ over $E _{v}$, for each $v \in M(E)$. This occurs
if and only if $c$ is presentable by $f$ over $E _{v}$, for each $v$
lying in the union $\Delta (R/E)$ of the sets $S _{p} (E)$, $p \in
P(M/E)$.
\end{coro}
\par
\medskip
\begin{proof} The assertion that $c \in N(R/E)$ is obviously
equivalent to the one that $c$ is presentable by $f$ over $F$. Since
$N(M/E) \subseteq N(R/E)$ (cf. \cite{L1}, Ch. VIII, Sect. 5), Lemmas
5.2 and 5.3 indicate that $N(R/E)$ includes as a subgroup the set
$\Sigma _{\varepsilon } = \{\lambda \in E ^{\ast }\colon  v(\lambda
- 1) < \varepsilon $, $v \in \Delta (R/E)\}$, for some $\varepsilon
> 0$ depending on $M/E$. Suppose now that $c$ is presentable by $f$
over $E _{v}$, for each $v \in \Delta (R/E)$, and fix a positive
number $\varepsilon ^{\prime }$. The polynomial $f$ maps $E _{v}$
continuously into itself, for each $v \in \Delta (R/E)$, so the
approximation theorem implies the existence of an element $c
_{\varepsilon '} \in N(R/E)$ such that $v(c - c _{\varepsilon '}) <
\varepsilon ^{\prime }$, $v \in \Delta (R/E)$. It is therefore clear
that $c.c _{\varepsilon '} ^{-1} \in N(R/E)$ when $\varepsilon
^{\prime }$ is sufficiently small, so Corollary 5.4 is proved.
\end{proof}
\par
\medskip
\begin{coro} In the setting of Corollary 5.4, $N _{\rm loc}(R/E) \le
N(R/E)$.
\end{coro}

\par
\medskip
\begin{proof} Our argument relies on the known fact $N
_{E} ^{R} (\rho ) = \prod _{j=1} ^{t} N _{E _{v}} ^{R _{\omega
_{j}}} (\rho )$, for each $\rho \in R ^{\ast }$, $v \in M(E)$, where
$\omega _{j}\colon \ j = 1,..., s$, are the prolongations of $v$ on
$R$ (cf. \cite{L1}, Ch. XII, Proposition 10). Also, it follows from
the approximation theorem that if $c _{j} \in E ^{\ast } \cap N(R
_{\omega _{j}}/E _{v})$: $j = 1,..., t$, and $\varepsilon $ is a
real positive number, then there is $\rho _{v} \in R ^{\ast }$, such
that $v(c _{j} - N _{E _{v}} ^{R _{\omega _{j}}} (\rho _{v})) <
\varepsilon $, for each index $j$. This implies that if $c \in N
_{\rm loc} (R/E)$, then one can find an element $\alpha
_{\varepsilon } \in R ^{\ast }$ so that $v(N _{E} ^{R} (\alpha
_{\varepsilon }) - c) < \varepsilon $, for every $v \in \Delta
(R/E)$. Therefore, Lemmas 5.2 and 5.3 yield $N _{E} ^{R} (\alpha
_{\varepsilon }).c ^{-1} \in N(M/E)$, in case $\varepsilon $ is
sufficiently small. Hence, $N _{\rm loc} (R/E) \subseteq N(R/E)$, as
claimed.
\end{proof}
\par
\medskip
To prove Theorem 2.3 we also need the following lemma.
\par
\medskip
\begin{lemm} Let $E$ be an algebraic extension of a global
field $E _{0}$ and $R$ a finite separable extension of $E _{v}$, for
some $v \in M(E)$. Also, let $E ^{\prime }$ be the subfield of $E
_{v}$ defined in (1.7), $R ^{\prime }$ the separable closure of $E
^{\prime }$ in $R$, and $R _{0}$ the maximal abelian extension of $E
_{v}$ in $R$. Then $N(R/E _{v}) = N(R _{0}/E _{v})$ and $N(R
^{\prime }/E ^{\prime }) = E ^{\prime \ast } \cap N(R/E _{v})$.
\end{lemm}

\par
\medskip
\begin{proof} It follows from (1.7) and \cite{Ch2}, Theorem~8.1,
that $E ^{\prime }$ is separably closed in $E _{v}$. Therefore, for
each finite extension $L$ of $E _{v}$ in $E _{v,{\rm sep}}$, the
norm map $N _{E _{v}}^{L}$ is a prolongation of $N _{E'} ^{L'}$, $L
^{\prime }$ being the separable closure of $E ^{\prime }$ in $L$.
Considering $R$ with its unique absolute value extending the
continuous prolongation $\bar v$ of $v$ on $E _{v}$, one obtains as
in the proof of Lemma 5.3 (ii) that $N(R/E _{v})$ includes the set
$\{\lambda \in E _{v} ^{\ast }\colon \bar v(\lambda - 1) <
\varepsilon \}$, for some real number $\varepsilon
> 0$. These observations show that $E _{v} ^{\ast } = E ^{\prime
\ast }N(R/E _{v})$ and $N(R ^{\prime }/E ^{\prime }) = E ^{\prime
\ast } \cap N(R/E _{v})$. Let $\Sigma $ be the normal closure of $R$
in $E _{v,{\rm sep}}$ over $E _{v}$, and let $F ^{\prime }$ be the
algebraic closure of $E ^{\prime }$ in $F$, for each $F \in I(\Sigma
/E _{v})$. It is clear from (1.3) and (1.4) that $\Sigma /E _{v}$ is
a Galois extension and the mapping of $I(\Sigma /E _{v})$ into
$I(\Sigma ^{\prime }/E ^{\prime })$ by the rule $F \to F ^{\prime }$
is bijective. Applying (1.4), (1.7) (ii) and Galois theory to
$\Sigma ^{\prime }/E ^{\prime }$, one also obtains that $R _{0}
^{\prime }$ is the maximal abelian extension of $E ^{\prime }$ in $R
^{\prime }$. Since $E ^{\prime }$ is quasilocal and $\rho _{E'/L'}$
is surjective, for every finite extension $L ^{\prime }$ of $E
^{\prime }$ (the latter follows from (1.5) (iv) and the
Albert-Hochschild theorem, see \cite{S}, Ch. II, 2.2), this enables
one to deduce from \cite{Ch6}, Theorem~1.1 (i), that $N(R ^{\prime
}/E ^{\prime }) = N(R _{0} ^{\prime }/E ^{\prime })$. Observing now
that $N(R _{0}/E _{v}) = N(R _{0} ^{\prime }/E ^{\prime })N(R/E
_{v}) = N(R ^{\prime }/E ^{\prime })N(R/E _{v})$, one concludes that
$N(R _{0}/E _{v}) \subseteq N(R/E _{v})$. The inclusion $N(R/E _{v})
\subseteq N(R _{0}/E _{v})$ follows at once from the transitivity of
norm maps in towers of finite separable extensions (cf. \cite{L1},
Ch. VIII, Sect. 5), so Lemma 5.6 is proved.
\end{proof}
\par
\medskip
{\it Proof of Theorem 2.3 and Proposition 2.4.} In view of Theorem
2.1 and Lemma 3.1, $E$ satisfies the conditions of Theorem 5.1. This
reduces Proposition 2.4 to a special case of Corollary 5.4, and
enables one to deduce from Corollary 5.5 that $N _{\rm loc} (R/E)
\subseteq N(R/E)$. It follows from Lemma 3.3 and the former part of
(1.5) (iii) that $E _{v(p)}$ admits local $p$-class field theory,
for every $p \in P(E)$. Since $N _{E} ^{R} (\rho ) = \prod _{v(p)'
\in \Sigma _{p}} N _{E _{v(p)}} ^{R _{v(p)'}} (\rho )$, for each
$\rho \in R ^{\ast }$, these observations, combined with Theorem 2.1
and Lemma 5.6, prove that $N(R/E) = N(F _{1}/E)$ and $N _{\rm loc}
(R/E) = N(F _{2}/E)$, where $F _{j} \in I(\overline E _{0}/E)$, $j =
1, 2$, are determined by Theorem 2.3 (iii) and (iv). They also show
that $F _{1} \subseteq F _{2}$ and $P(F _{2}/E) \subseteq P(M/E)$.
Now Theorem 2.3 (ii) can be deduced from Burnside-Wielandt's
characterization of nilpotent finite groups (cf. \cite{KM}, Ch. 6,
Sect. 2), Galois theory and Lemma 3.1 (i). As the latter part of
Theorem 2.3 (i) is contained in Remark 3.6, Theorem 2.3 is proved.
\par
\medskip
We conclude this Section with an example of a field extension
$R/E$ satisfying the conditions of Theorem 2.3, for which $N _{\rm
loc} (R/E) \neq N(R/E)$.
\par
\medskip
\begin{exam} Let $w(p)$ the $p$-adic absolute value of $\QQ $, for
each $p \in \PP $, $W = \{w(p)\colon \ p \in \PP \}$, $E$ a minimal
element of $\Omega (\QQ , \PP , W)$, $V = \{v(p)\colon \ p \in \PP
\}$ a characteristic system of $E$, and $R = E(\alpha )$, where
$\alpha \in \overline {\QQ }$ is a root of the polynomial $f(X) = X
^{4} + 8X ^{2} + 256X - 20$. Denote by $g(X) = X ^{3} - 16X ^{2} +
144X + 256 ^{2}$ the resolution polynomial of $f(X)$, and put
$\overline {\QQ } _{2} = \QQ _{2,{\rm sep}}$. It is easily verified
that $f(X)$ has roots in $\QQ _{3}$ and in the subfields $\QQ _{2}
(\sqrt{2})$ and $\QQ _{2} (\sqrt{-10})$ of $\overline {\QQ } _{2}$,
but has no roots in $\QQ _{2}$ (cf. \cite{L2}, Ch. II, Proposition
2). It is similarly obtained that $g(X)$ is irreducible over $\QQ
_{3}$ and the discriminant $d(g) = d(f)$ lies in $\QQ _{3} ^{\ast
2}$. Observing also that $-5d(f) \in \QQ _{2} ^{\ast 2}$, whereas
$-5 \not\in \QQ _{2} ^{\ast 2}$, one proves that $f(X)$ and $g(X)$
are irreducible over $\QQ $ with Galois groups $\mathcal{G} _{f}
\cong {\rm Sym} _{4}$ and $\mathcal{G} _{g} \cong {\rm Sym} _{3}$.
These results, combined with Theorem 2.2 (iii), imply that $g(X)$
and $f(X)$ remain irreducible over $E$ with the same Galois groups.
They also show that $v(p)$ has two nonequivalent prolongations $v
_{1} (p)$ and $v _{2} (p)$ on $R$, for $p = 2, 3$. Hence, by Theorem
2.3, $N(R/E) = E ^{\ast }$ and $N _{\rm loc} (R/E)$ is a subgroup of
$E ^{\ast }$, such that $E ^{\ast }/N _{\rm loc} (R/E)$ is noncyclic
of order $12$. Moreover, it becomes clear that an integer $\beta $
not divisible by $3$ lies in $N _{\rm loc} (R/E)$ if and only if
$\beta $ or $-2\beta  \in \QQ _{2} ^{\ast 2}$. In particular, $-20
\not\in N _{\rm loc} (R/E)$.
\end{exam}

\end{document}